\documentclass[11pt, a4paper,reqno]{amsart}
\usepackage{color} \definecolor{bleu_sombre}{rgb}{0,0,0.6}  \definecolor{rouge_sombre}{rgb}{0.8,0,0}\definecolor{vert_sombre}{rgb}{0,0.6,0}
\usepackage[plainpages=false,colorlinks,linkcolor=bleu_sombre,citecolor=rouge_sombre,urlcolor=vert_sombre,breaklinks]{hyperref}
\usepackage[english]{babel}
\usepackage{amsmath,amssymb,amsthm,graphicx,amsfonts,url,color,enumerate,dsfont,stmaryrd,mathabx}

\theoremstyle{plain}
\newtheorem{theorem}{{Theorem}}[section] 
\newtheorem*{theorem*}{{Theorem}}
\newtheorem{proposition}[theorem]{Proposition}
\newtheorem*{proposition*}{Proposition}
\newtheorem{corollary}[theorem]{Corollary}
\newtheorem*{corollary*}{Corollary}
\newtheorem{lemma}[theorem]{Lemma}
\newtheorem*{lemma*}{Lemma}

\theoremstyle{definition}
\newtheorem{definition}[theorem]{Definition}
\newtheorem*{definition*}{Definition}
\theoremstyle{remark}
\newtheorem{remark}[theorem]{Remark}

\makeatletter

\@addtoreset{equation}{section}  
\makeatother

\newcommand {\limt}[2]{\xrightarrow[#1 \to #2]{}}

\newcommand{\abs}[1]{\left\vert #1\right\vert}        
\newcommand{\nr}[1]{\left\Vert #1\right\Vert}         

\newcommand{\innp}[2]{\left\langle #1 , #2 \right \rangle}        
\newcommand{\pppg}[1] {\left< #1 \right>}
\newcommand{\set}[1]{\left\{ #1 \right\}}
\newcommand{\Ii}[2] {\left\{ #1,\dots,#2 \right\}}
\renewcommand{\leq}{\leqslant}	\renewcommand{\geq}{\geqslant}

\newcommand{\inv}{^{-1}}

\newcommand{\1}{\mathds 1}
\newcommand{\st}{\,:\,}

\renewcommand{\Re}{\mathsf{Re}}
\renewcommand{\Im}{\mathsf{Im}}

\newcommand{\Dom}{\mathsf{Dom}}

\newcommand{\Ran}{\mathsf{Ran}}
\newcommand{\Id}{\mathsf{Id}}

\newcommand{\divg}{\mathop{\rm{div}}\nolimits}

\newcommand{\R}{\mathbb{R}}		\newcommand{\C}{\mathbb{C}}
\newcommand{\N}{\mathbb{N}}		
\newcommand{\DD}{\mathbb{D}}

\renewcommand{\a}{\alpha}\renewcommand{\b}{\beta}\newcommand{\g}{\gamma}\renewcommand{\d}{\delta}\newcommand{\D}{\Delta}\newcommand{\e}{\varepsilon}\newcommand{\z}{\zeta} \newcommand{\y}{\eta}\renewcommand{\th}{\theta}\newcommand{\Th}{\Theta}\renewcommand{\k}{\kappa}\renewcommand{\l}{\lambda}\newcommand{\m}{\mu}\newcommand{\x}{\xi}\newcommand{\s}{\sigma}\newcommand{\f}{\varphi}\newcommand{\vf}{\phi}

\newcommand{\Cc}{{\mathcal C}}\newcommand{\Dc}{{\mathcal D}}\newcommand{\Hc}{{\mathcal H}}\newcommand{\Ic}{{\mathcal I}}\newcommand{\Kc}{{\mathcal K}}\newcommand{\Lc}{{\mathcal L}}\newcommand{\Nc}{{\mathcal N}}\newcommand{\Rc}{{\mathcal R}}\newcommand{\Sc}{{\mathcal S}}\newcommand{\Vc}{{\mathcal V}}\newcommand{\Wc}{{\mathcal W}}

\newcommand{\stepp}{\noindent {\bf $\bullet$}\quad }

\newcommand{\detail}[1]
{
}

\newcommand{\diff}{\, \mathrm d}

\usepackage{mathrsfs}

\begin{document}

\newcommand{\HH}{\mathscr H}\newcommand{\LL}{\mathscr L}\newcommand{\EE}{\mathscr E} \newcommand{\SSS}{\mathscr S}

\newcommand{\Pa}{P_{\mathsf {a}}}
\newcommand{\Ra}{R_{\mathsf {a}}}
\newcommand{\Rs}{R_{\mathsf {s}}}
\newcommand{\PRinf}{P_{\mathsf{R},\infty}}

\newcommand{\Opw}{{\mathop{\rm{Op}}}_h^w}
\newcommand{\Opwn}{{\mathop{\rm{Op}}}_{h_n}^w}

\newcommand{\nul}{{\nu_l}}
\newcommand{\nur}{{\nu_r}}
\newcommand{\nus}{{\nu_*}}
\newcommand{\anul}{{\abs \nul}}
\newcommand{\anur}{{\abs \nur}}
\newcommand{\anus}{{\abs \nus}}
\newcommand{\tnul}{{\tilde \nu_l}}
\newcommand{\tnur}{{\tilde \nu_r}}
\newcommand{\tnus}{{\tilde \nu_*}}
\newcommand{\mul}{{\mu_l}}
\newcommand{\mur}{{\mu_r}}
\newcommand{\Rci}{\Rc^\y}
\newcommand{\tdelta}{\s}
\newcommand{\Phil}{\Phi_l(z)}
\newcommand{\Phir}{\Phi_r(z)}
\newcommand{\tPhil}{\tilde \Phi_l(z)}
\newcommand{\tPhir}{\tilde \Phi_r(z)}
\newcommand{\tVc}{\tilde \Vc}

\newcommand{\ww}{w}
\newcommand{\ad}{\mathsf{ad}}\newcommand{\Ad}{\mathsf{Ad}}

\newcommand{\aright}{\underrightarrow{\mathfrak a}}
\newcommand{\aleft}{\underleftarrow{\mathfrak a}}

\newcommand{\tWc}{\widetilde \Wc}

\newcommand{\DI}{\mathbb D_{\mathsf I}}
\newcommand{\DR}{\mathbb D_{\mathsf R}}

\newcommand{\Qp}{Q_+} \newcommand{\Qpp}{Q_\bot^+}

\newcommand{\PR}{P_{\mathsf{R}}}

\newcommand{\jj}{\mathfrak j}
\renewcommand{\ll}{\mathfrak l}

\newcommand{\dimm}[1]{
} 

\newcommand{\h}{h}

\newcommand{\A}{A} \newcommand{\Ao}{{A_0}}
\newcommand{\Asf}{\mathsf A}
\author{R. Fahs  }
\author{J. Royer  }
\address{Institut de Math\'ematiques de Toulouse -- UMR 5219 -- Universit\'e de Toulouse, F-31062 Toulouse Cedex 9, France.}
\email{Rayan.Fahs@math.univ-toulouse.fr}
\email{julien.royer@math.univ-toulouse.fr}
\title[Local decay for the damped wave equation]{Local decay and asymptotic profile for the damped wave equation in the asymptotically Euclidean setting}

\begin{abstract}
We prove local decay estimates for the wave equation in the asymptotically Euclidean setting. In even dimensions we go beyond the optimal decay by providing the large time asymptotic profile, given by a solution of the free wave equation. In odd dimensions, we improve the best known estimates. In particular, we get a decay rate that is better than what would be the optimal decay in even dimensions.
The analysis mainly relies on a comparison of the corresponding resolvent with the resolvent of the free problem for low frequencies. Moreover, all the results hold for the damped wave equation with short range absorption index.
\end{abstract}

\maketitle

\section{Introduction}
We consider on $\R^d$, $d\geq 3$, the (possibly damped) wave equation
\begin{equation} \label{wave}
\begin{cases}
 \partial_t^2 u + P u +a(x)\partial_t u  = 0, & \quad \text{on } \R_+ \times \R^d,\\
(u,\partial_t u)_{|t = 0} = (f,g), & \quad \text{on } \R^d,
\end{cases}
\end{equation}
where $(f,g) \in H^1\times L^2$, the operator $P$ is a general Laplace operator on $\R^d$, close to the free Laplacian at infinity, and the absorption index $a(x)$ is small at infinity.

More precisely, $P$ is of the form
\begin{equation} \label{def-L}
P = - \frac 1 {w(x)} \divg G(x) \nabla,
\end{equation}
where the density $w(x)$ and the symmetric matrix $G(x)$ are smooth and uniformly positive functions: there exist $C_G,C_w \geq 1$ such that, for all $x \in \R^d$ and $\xi \in \R^d$,
\begin{equation} \label{eq:CGCw}
C_G \inv \abs \xi^2 \leq \innp{G(x) \xi}{\xi}_{\R^d} \leq C_G \abs \xi^2 \quad \text{and} \quad C_w\inv \leq w(x) \leq C_w.
\end{equation}
We assume that $P$ is associated to a long range perturbation of the flat metric. This means that $G(x)$ and $w(x)$ are long range perturbations of $\Id$ and 1, respectively, in the sense that for some $\rho_0 \in ]0,1]$ there exist constants $C_\a > 0$, $\a \in \N^d$, such that for all $x \in \R^d$,
\begin{equation} \label{hyp-swa}
\big| \partial^\a (G(x) - \Id) \big| + \big|\partial^\a (w(x) - 1)\big| \leq C_\a \pppg x^{-\rho_0 - \abs \a}.
\end{equation}
Here and everywhere below, we use the standard notation $\pppg x = (1 + \abs x^2)^{\frac 12}$. We also denote by $\D_G$ the Laplace operator in divergence form corresponding to $G$:
\[
\D_G = \divg G(x) \nabla.
\]

This definition of $P$ includes in particular the cases of the free Laplacian, a Laplacian in divergence form, or a Laplace-Beltrami operator. We recall that the Laplace-Beltrami operator associated to a metric $\mathsf g = (g_{j,k})_{1\leq j,k\leq d}$ is given by
\[
P_{\mathsf g}  = - \frac 1 {\abs{g(x)}^{\frac 12}} \sum_{j,k=1}^d \frac {\partial}{\partial x_j} \abs{g(x)}^{\frac 12} g^{j,k}(x) \frac {\partial}{\partial x_k},
\]
where $\abs{g(x)} = \abs{\det(\mathsf g(x))}$ and $(g^{j,k}(x))_{1\leq j,k\leq d} = \mathsf g(x)\inv$. Then $P_{\mathsf g}$ is of the form \eqref{def-L} with $w = \abs{g}^{\frac 12}$ and $G = \abs{g}^{\frac 12} \mathsf g\inv$.\\

On the other hand, the absorption index $a(x)$ is smooth, bounded, takes non-negative values, and is of short range: there exist $C_\a > 0$, $\a \in \N^d$, such that, for all $x\in \R^d$,
\begin{equation} \label{hyp-a}
\big| \partial^\a a(x) \big|\leq C_\a \pppg x^{-1-\rho_0 - \abs \a}.
\end{equation}
In particular, $a$ can be identically 0, so our setting includes the undamped wave equation.\\

The solution of the wave equation is explicit for the free case
\begin{equation} \label{wave-free}
\begin{cases}
   \partial_t^2 u_0 -\D u_0= 0, & \quad \text{on } \R_+ \times \R^d,\\
(u_0,\partial_t u_0)_{|t = 0} = (f_0,g_0), & \quad \text{on } \R^d,
\end{cases}
\end{equation}
see for instance \cite{evans, CourantHilbert}. In particular, if the dimension $d$ is odd, the wave propagates at speed 1 (this is the strong Huygens Principle). Then, given $R > 0$ and $(f_0,g_0)$ supported in the ball $B(R)$ of radius $R$, the solution $u_0$ of \eqref{wave-free} vanishes on the ball $B(R)$ for any time $t \geq 2R$. The situation is different in even dimension, but the wave still escapes to infinity. More precisely, we have the estimate
\begin{equation} \label{eq:loc-dec-intro}
\nr{u_0(t)}_{L^2(B(R))} \lesssim \pppg t^{-d} \nr{f_0}_{L^2(B(R))}+ \pppg t^{1-d} \nr{g_0}_{L^2(B(R))}.
\end{equation}
We can separate the contribution of $f_0$ and $g_0$ by writing
\begin{equation} \label{eq:u0}
u_0(t) = \cos(t\sqrt{-\D}) f_0 + \frac {\sin(t \sqrt{-\D})}{\sqrt{-\D}} g_0,
\end{equation}
and then we have
\begin{equation} \label{eq:cos}
\nr{\cos(t\sqrt{-\D}) f_0}_{L^2(B(R))} \lesssim \pppg t^{-d} \nr{f_0}_{L^2(B(R))},
\end{equation}
and
\begin{equation} \label{eq:sin}
\nr{\frac {\sin(t \sqrt{-\D})}{\sqrt{-\D}} g_0}_{L^2(B(R))} \lesssim \pppg t^{1-d} \nr{g_0}_{L^2(B(R))}.
\end{equation}
Moreover, as observed in \cite{BoucletBur21}, these estimates are optimal if the integral of $f_0$ (respectively $g_0$) is not 0.

Finally, we recall that in dimension 1, the solution of the free wave equation is given by the d'Alembert formula
\[
u_0(t,x) = \frac {f_0(x+t) + f_0(x-t)}{2} + \frac 12 \int_{x-t}^{x+t} g_0(s) \diff s.
\]
We easily see that if $f_0$ and $g_0$ are supported in $]-R,R[$ then for $t \geq 2R$, the first term vanishes in $]-R,R[$ (as in any odd dimension), while the second term is equal to a constant (half of the integral of $g_0$).\\

Our purpose in this paper is to prove such local decay estimates for the solution of the perturbed wave equation \eqref{wave}. 
\\

We will work from a spectral point of view, and study this time-dependent problem with a frequency-dependent analysis. In particular, it is well known that the contributions of high and low frequencies play very different roles.

A high frequency, a wave propagates along the classical trajectories. In our setting, we set on $\R^{2d} \simeq T^* \R^d$
\[
p(x,\xi) = w(x)\inv \innp{G(x) \xi}{\xi}_{\R^d}.
\]
The classical rays of light are the solutions of the corresponding Hamiltonian problem
\begin{equation*}
\begin{cases}
x'(t;x_0,\xi_0) = \partial_\xi p \big(x(t;x_0,\xi_0),\xi(t;x_0,\xi_0)\big),\\
\xi'(t;x_0,\xi_0) = - \partial_x p \big(x(t;x_0,\xi_0),\xi(t;x_0,\xi_0) \big),\\
(x,\xi)(0;x_0,\xi_0) = (x_0, \xi_0).
\end{cases}
\end{equation*}
For example, if $P$ is the Laplace-Beltrami operator $P_{\mathsf g}$, then the rays of light are the geodesics of the metric $\mathsf g$. The geometry of these rays of light plays an important role for the analysis of high frequencies.

The behavior of the contribution of low frequencies is completely different. Geometry does not play any particular role, but the setting at infinity is crucial. In particular, the fact that our problem looks like the free wave equation at infinity (in the sense given by Assumptions \eqref{hyp-swa} and \eqref{hyp-a}) will be important. This is even more important given that the rate for local decay is limited by the contribution of low frequencies.\\

This question of local decay has a long history for the wave and Schr\"odinger equations.

The first results are about the undamped wave equation. We refer for instance to \cite{morawetz61}, where a multiplier method is used for the free wave outside a star-shaped obstacle in dimension 3. Then exponential decay is proved in \cite{laxmp63} via an analysis of the corresponding semigroup. See also \cite{LaxPhi62,LaxPhi72,LaxPhillips89}.

We have said that the contribution of high frequencies follows the classical rays of light. Without damping, it is then natural that waves escape to infinity if these rays of light do. The following non-trapping condition is then very important for the local decay of the undamped wave equation:
\begin{equation} \label{eq:non-trapping}
\forall (x_0,\xi_0) \in p\inv(\{1\}), \quad \abs {x(t;x_0,\xi_0)} \limt t {\pm \infty} \infty.
\end{equation}
It is proved in \cite{ralston69} that this non-trapping condition is necessary to have uniform local energy decay. Local energy decay outside non-trapping obstacles is considered in \cite{morawetz75,strauss75,morawetzrs77}.

In \cite{Vainberg75}, the properties of the time-dependent problem are deduced from the analysis of the stationary problem. Another important step is the analysis of  \cite{melrose79}, based on the propagation of singularities of \cite{melroses78}. See also \cite{Kawashita93,Vodev99}. We also refer to \cite{burq98} for the logarithmic decay (with loss of regularity) outside any (in particular, trapping) compact smooth obstacle.

The recent papers deal simultaneously with the local decay for the Schr\"odinger and wave equations, which are now similar from this spectral point of view. We refer to \cite{Bouclet11,BonyHaf12} for estimates with an $\e$-loss on an asymptotically Euclidean setting. The $\e$-loss has finally been removed in \cite{BoucletBur21}. The method does not see the parity of the dimension, so this final result is optimal for Schr\"odinger or for the wave in even dimension, but not for the wave in odd dimension. However, a better result is obtained in \cite{BonyHaf13} when the metric goes faster to the flat metric at infinity.
In these works, the main contribution is the analysis of low frequency resolvent estimates (see also \cite{wang06,DerezinskiSki09,Bouclet11b,bonyh10b}).

High frequency resolvent estimates were already understood for the Schr\"odinger operator in close settings. See, for instance, \cite{robertt87,Robert92,Burq02}.\\

Here, we also consider the damped wave equation. Stabilization of the wave equation also has a long history on compact domains. In particular, it is known that the global energy decays uniformly (hence exponentially) under the Geometric Control Conditions (all the classical trajectories go through the damping region), and we get weaker results with loss of regularity when this condition is not satisfied. See, for instance, \cite{RauchTay74,BardosLebRau92,Lebeau96,LebeauRob97,BurqHit07,christianson07,leautaudl,BurqGer20}.\\

In unbounded domains we have additional difficulties, in particular due to the contribution of low frequencies.

For the damped wave equation on an unbounded domain, the size of the solution on a compact now has two reasons to decay. Either because it escapes to infinity, or because it is dissipated. The expected corresponding condition on classical trajectories is that they should all escape to infinity or go through the damping region. This means that we can allow trapped trajectories if they are damped. We set
\[
\Omega_{\mathrm b} = \set{(x_0,\x_0) \in p\inv(\set 1) \, : \, \sup_{t \in \R} \abs {x(t;x_0,\xi_0)} < +\infty}.
\]
Then the condition on classical trajectories reads
\begin{equation} \label{hyp:damping}
\forall (x_0,\xi_0) \in \Omega_{\mathrm b} , \exists t \in \R, \quad a(x(t,x_0,\xi_0)) > 0.
\end{equation}
We refer to \cite{AlouiKhe02} for the damped wave equation in odd dimensions, in an exterior domain, and with a compactly supported damping, and to \cite{khenissi03} for even dimensions. The damped wave equation in the asymptotically Euclidean setting has been studied in \cite{royer-mourre,BoucletRoy14,royer:dld-energy-space}, where the local energy decay with $\e$-loss has been obtained.\\

The present work improves the results of \cite{BoucletRoy14,royer:dld-energy-space} about the damped wave equation, and also the sharp estimates of \cite{BoucletBur21} for the undamped case.
More precisely, in a general setting including the damped case, we go beyond the optimal estimate of \cite{BoucletBur21}. For this, as is done in \cite{Royer24} for the Schr\"odinger equation, we estimate the difference between the solution $u(t)$ of \eqref{wave} and a solution $u_0(t)$ of the free wave equation \eqref{wave-free}. We prove that this difference decays faster than the rates given by \eqref{eq:cos}-\eqref{eq:sin}.

To state the main result of this paper, we introduce some notation. Instead of considering compactly supported solutions and estimating the solution in a compact, we work in weighted spaces. For $\d \in \R$, we set $L^{2,\d} = L^2(\pppg x^{2\d} \diff x)$ and we denote by $H^{1,\d}$ the corresponding Sobolev space, with norm defined by
\[
\nr{u}_{H^{1,\d}}^2 = \nr{\nabla u}_{L^{2,\d}}^2 + \nr{u}_{L^{2,\d}}^2.
\]
Our results rely on the damping condition \eqref{hyp:damping}, saying that all the geodesics go through the damping region or go to infinity. We recall that $\rho_0 \in ]0,1]$ measures the decay to the free setting in \eqref{hyp-swa} and \eqref{hyp-a}.

\begin{theorem} \label{th:loc-decay}
Assume that the damping condition \eqref{hyp:damping} holds. Let $\rho_1 \in ]0,\rho_0[$ and $\d > d + \frac 52$. There exists $C > 0$ such that for $f \in H^{1,\d}$, $g \in L^{2,d}$ and $t \geq 0$, we have
\[
\nr{u(t) - u_0(t)}_{L^{2,-\d}} \leq C \pppg t^{-d-\rho_1} \nr{f}_{H^{1,\d}} +  C \pppg t^{1-d-\rho_1} \nr{af + g}_{L^{2,\d}},
\]
where $u$ is the solution of the damped wave equation \eqref{wave} and $u_0$ is the solution \eqref{eq:u0} of the free wave equation \eqref{wave-free} with initial condition 
\[
(f_0,g_0) = (wf,awf +wg).
\]
\end{theorem}

With this result, we generalize in particular the optimal estimates of \cite{BoucletBur21} to the undamped wave equation.

\begin{theorem} \label{th:loc-dec-even}
Assume that $d$ is even and that the damping condition \eqref{hyp:damping} holds. Let $R > 0$. Then there exists $C > 0$ such that for $f \in H^{1}$ and $g \in L^{2}$ supported in $B(R)$ we have, for all $t \geq 0$, 
\[
\nr{u(t)}_{L^{2}(B(R))} \leq C \pppg t^{-d} \nr{f}_{H^1} +  C \pppg t^{1-d} \nr{af + g}_{L^2},
\]
where $u$ is the solution of the damped wave equation \eqref{wave}. Moreover, the estimate is optimal if $d$ is even and $\int_{\R^d} f \diff x \neq 0$ or $\int_{\R^d} (af+g) \diff x \neq 0$.
\end{theorem}

The conclusions of Theorem \ref{th:loc-decay} are different in even and odd dimensions. In even dimension, it proves that $u(t)$ is equal to $u_0(t)$ up to a smaller rest. This not only means that $u(t)$ decays exactly like $u_0(t)$ (as stated in Theorem \ref{th:loc-dec-even}, which is already known in the undamped case but new with damping) but also gives the leading term for the large time asymptotic expansion of $u(t)$ in weighted spaces: $u(t)$ looks like $u_0(t)$ on compacts for large time.

In odd dimension, the situation is different, since the solution of the free problem decays very fast. In this case, we do not get the asymptotic profile of the solution $u(t)$ but we improve the local decay. Even if we do not get an optimal result, going beyond the $t^{1-d} / t^{-d}$ rate of decay is a very important improvement. Indeed, the best results about local decay are based on the Mourre commutators method (see Section \ref{sec:Mourre} below), which does not see the parity of the dimension. Thus, even in odd dimension, it cannot give a result better than what is the optimal decay in even dimension. Here, we will use the same commutators method, but by comparing the solution with the solution of the free wave equation, we reintroduce a difference between odd and even dimensions. This gives the following result for odd dimensions.

\begin{theorem} \label{th:loc-dec-odd}
Assume that $d$ is odd and that the damping condition \eqref{hyp:damping} holds. Let $\rho_1 \in ]0,\rho_0[$ and $R > 0$. There exists $C > 0$ such that for $f \in H^{1}$ and $g \in L^{2}$ supported in $B(R)$ we have, for all $t \geq 0$,
\[
\nr{u(t)}_{L^{2}(B(R))} \leq C \pppg t^{-d-\rho_1} \nr{f}_{H^1} +  C \pppg t^{1-d-\rho_1} \nr{g}_{L^{2}},
\]
where $u$ is the solution of the damped wave equation \eqref{wave}.
\end{theorem}

Notice, however, that for the undamped wave equation, if $G - \Id$ decays fast enough at infinity, it is proved in \cite{BonyHaf13} that the time decay for the local energy improves with the spatial decay rate toward the free metric. More precisely, for any $\rho_0,\e > 0$, the local energy decays like $t^{-\rho_0}$ if the metric converges like $\abs{x}^{-\rho_0-2-\e}$ toward the Euclidean metric. We use a much weaker assumption in this paper.
\\

Theorem \ref{th:loc-decay} will be proved from a spectral point of view. Given $\mu > 0$, we can write
\begin{equation} \label{eq:u-Rz}
u(t) = \frac 1 {2\pi} \int_{\Im(z) = \mu} e^{-itz} R(z) ( awf - izwf + wg ) \diff z,
\end{equation}
where for $\Im(z) > 0$, we have set
\begin{equation} \label{def:R}
R(z) = \big(-\D_G - izaw- z^2 w\big)\inv.
\end{equation}
Due to the (growing) exponential factor in \eqref{eq:u-Rz}, our purpose is to prove estimates for the right-hand side that are uniform in $\mu > 0$, and to let $\mu$ go to 0. Moreover, since we want to prove time decay for $u(t)$, we have to estimate the derivatives of the integrand in the right-hand side of \eqref{eq:u-Rz}.

This implies that we have to prove estimates for $R(z)$ and its derivatives near the real axis in weighted spaces (limiting absorption principle). See details in Section \ref{sec:time-decay}. The main difficulties are due to the contribution of high frequencies (for $z$ large) and low frequencies (for $z$ near 0).\\

For high frequencies, we have the following estimates (see \cite[Th. 1.5]{BoucletRoy14}).

\begin{theorem} \label{th:high-freq}
Assume that the damping condition \eqref{hyp:damping} holds. Let $n \in \N$ and $\d > n + \frac 12$. Let $\tau_0 > 0$. There exists $c > 0$ such that for $z \in \C$ with $\abs{\Re(z)} \geq \tau_0$ and $\Im(z) > 0$, we have
\[
\nr{\pppg x^{-\d} R^{(n)}(z) \pppg x^{-\d}}_{\Lc(L^2)} \leq \frac c {\abs z}.
\]
\end{theorem}

With Theorem \ref{th:high-freq}, we will check that the contribution of high frequencies decays very fast. As said above, the rate of decay for the local decay is governed by the contribution of low frequencies. This is the main issue of this paper.

Since our purpose is to compare the solution of the perturbed wave equation with a solution of the free problem, we compare $R(z)$ with
\begin{equation} \label{def:R0}
R_0(z) = (-\D - z^2)\inv.
\end{equation}
The main part of this paper will be the proof of the following low frequency estimates.

\begin{theorem} \label{th:low-freq}
Let $\rho_1 \in ]0,\rho_0[$. Let $n \in \N$ and $\d > n + \frac 32$. There exists $C > 0$ such that for $z \in \C$ with $\abs z \leq 1$ and $\Im(z) > 0$, we have
\[
\nr{\pppg x^{-\d} \big(R^{(n)}(z)- R_0^{(n)}(z) \big) \pppg x^{-\d}}_{\Lc(L^{2})} \leq C \abs z^{\min(d+ \rho_1-n-2,0)}.
\]
\end{theorem}

As for the time decay problem, this result gives the optimal resolvent estimate in even dimensions and improves the best estimates in odd dimensions.

\begin{corollary}
Let $n \in \N$ and $\d > n + \frac 32$. There exists $C > 0$ such that for $z \in \C$ with $\abs z \leq 1$ and $\Im(z) > 0$, we have
\[
\nr{\pppg x^{-\d} R^{(n)}(z) \pppg x^{-\d}}_{\Lc(L^2)} \leq C \abs z^{\min(d-n-2,0)}.
\]
\end{corollary}

\begin{corollary}
Let $\rho_1 \in ]0,\rho_0[$. Assume that $d$ is odd. Let $n \in \N$ and $R > 0$. There exists $C > 0$ such that for $z \in \C$ with $\abs z \leq 1$ and $\Im(z) > 0$, we have
\[
\nr{\1_{B(R)} R^{(n)}(z) \1_{B(R)}}_{\Lc(L^2)} \leq C \abs z^{\min(d+\rho_1-n-2,0)}.
\]
\end{corollary}

\subsection*{Plan of the paper} The paper is organized as follows. In Section \ref{sec:time-decay}, we give some basic properties for the resolvent $R(z)$ and we derive the local decay stated in Theorem \ref{th:loc-decay} from the resolvent estimates of Theorems \ref{th:high-freq} and \ref{th:low-freq}. Then, it remains to prove Theorem \ref{th:low-freq}. In Section \ref{low-freq}, we outline the strategy of the proof, assuming some intermediate results. These results are proved in the remaining sections. In Section \ref{sec:hardy-ineq}, we show how the decay of the coefficients at infinity (see \eqref{hyp-swa} and \eqref{hyp-a}) and the weights $\pppg x^{-\d}$ provide smallness for low frequencies, and in Section \ref{sec:elliptic-regularity} we combine this observation with the elliptic regularity of $R(z)$ to prove some resolvent estimates. Finally, in Section \ref{sec:Mourre}, we recall the Mourre method and use it to prove estimates for $R(z)$ when $z$ is close to the real axis.

\section{From resolvent estimates to time decay} \label{sec:time-decay}

In this section, we check that the resolvent $R(z)$ (introduced in \eqref{def:R}) is well defined and we show how we can deduce Theorem \ref{th:loc-decay} from Theorems \ref{th:high-freq} and \ref{th:low-freq}. By density, it is enough to prove the estimates of $f$ and $g$ in the Schwartz space $\Sc$.\\

We set
\[
\C^\pm = \set{z \in \C \, : \, \pm\Re(z) > 0}, \quad \C_\pm = \set{z \in \C \, : \, \pm \Im(z) > 0}.
\]
We say that an operator $T$ on a Hilbert space $\Hc$ is accretive (resp. dissipative) if
\[
\forall \f \in \Dom(T), \quad \Re \innp{T\f}{\f}_\Hc \geq 0 \quad (\text{resp. } \Im \innp{T\f}{\f}_\Hc \leq 0).
\]
Moreover, $T$ is maximal accretive (resp. maximal dissipative) if some (hence any) $z \in \C^-$ (resp. $z \in \C_+$) belongs to the resolvent set of $T$.\\

For $z \in \C_+$ we consider on $L^2$ the operators
\begin{equation} \label{def:Pz}
P(z) = (-\D_G -izaw -z^2 w), \quad P_0(z)  = (-\D - z^2),
\end{equation}
with domain $\Dom(P(z)) = \Dom(P_0(z)) = H^2$. They can also be seen as bounded operators from $H^1$ to its dual $H\inv$.

\begin{proposition} \label{prop:Rz-base}
Let $z \in \C_+$.
\begin{enumerate}[\rm(i)]
\item The operator $P(z)$ is boundedly invertible and $R(z)$ --see \eqref{def:R}-- is a well defined bounded operator on $L^2$.
\item We have $R(z)^* = R(-\overline z)$.
\item The resolvent $R(z)$ extends to a bounded operator from $H^{s-1}$ to $H^{s+1}$ for any $s \in \R$.
\end{enumerate}
\end{proposition}

\begin{proof}
Let $\vartheta_z = \arg(-iz) \in \big]-\frac \pi 2,\frac \pi 2 \big[$. By \eqref{eq:CGCw}, there exists $c_0 > 0$ such that for all $u \in H^1$
\begin{eqnarray} \label{eq:Lax-Milgram}
\lefteqn{\Re \big < e^{-i\vartheta_z} P(z) u,u\big >}_{H^1,H\inv}\\
\nonumber
&& = \cos(\vartheta_z) \innp{G\nabla u}{\nabla u}_{L^2} + \abs z \innp{aw u}{u}_{L^2} + \cos(\vartheta_z) \abs z^2 \innp{wu}{u}_{L^2}\\
\nonumber
&& \geq c_0 \nr{u}_{H^1}^2.
\end{eqnarray}
By the Lax-Milgram Theorem, this implies that $e^{-it\vartheta_z} P(z)$ and hence $P(z)$ define boundedly invertible operators from $H^1$ to $H\inv$. Then $R(z)$ is well defined as an operator from $H\inv$ from $H^1$. By elliptic regularity, it can also be seen as a bounded operator from $L^2$ to $H^2$. More generally, by elliptic regularity, duality and interpolation, it defines a bounded operator from $H^{s-1}$ to $H^{s+1}$ for any $s \in \R$.
Finally, we have
\[
P(z)^* = -\D_G +i\bar z a w -\bar z^2 w = P(-\bar z),
\]
so $R(z)^* = R(-\bar z)$.
\end{proof}

Next, we recall that \eqref{wave} is well posed. We set $\HH = H^1 \times L^2$, and we define on $\HH$ the operator
\[
\Wc = \begin{pmatrix} 0 & w\inv \\ \D_G & -a \end{pmatrix},\quad \Dom(\Wc) = H^2 \times H^1.
\]

\begin{proposition} \label{prop:Wc}
The operator $\Wc$ generates a $C^0$-semigroup on $\HH$. Moreover, for $\e > 0$ there exists $M_\e > 0$ such that for all $t \geq 0$, we have
\[
\nr{e^{t\Wc}}_{\Lc(\HH)} \leq M_\e e^{t\e}.
\]
\end{proposition}

\begin{proof}
Let $\nu = {2\e / \|w^{-1/2}\|_{L^\infty}}$. We consider on $\HH$ the norm $\nr{\cdot}_{\HH,\nu}$ defined by
\[
\nr{(u,v)}_{\HH,\nu}^2 = \innp{G\nabla u} {\nabla u}_{L^2} + \nu^2 \nr{u}_{L^2}^2 + \innp{w\inv v}{v}_{L^2}, \quad (u,v) \in \HH.
\]
It is equivalent to the usual norm on $\HH$. For $U = (u,v) \in \Dom(\Wc)$, we have
\begin{align*}
\innp{\Wc U}{U}_{\HH,\nu}
& = \innp{G \nabla (w\inv v)}{\nabla u}_{L^2} + \nu^2 \innp{w\inv v}{u}_{L^2}\\
& + \innp{w\inv \D_G u}{v}_{L^2} - \innp{aw\inv v}{v}_{L^2},
\end{align*}
so
\begin{align*}
\Re \innp{\Wc U}{U}_{\HH,\nu}
& \leq \nu^2 \Re \innp{w\inv v}{u}_{L^2} \leq \frac {\nu \|w^{-1/2}\|_{L^\infty}}2 \big( \nu^2 \nr{u}_{L^2}^2 + \innp{w\inv v}{v}_{L^2} \big)\\
& \leq \e \nr{U}_{\HH,\nu}^2.
\end{align*}
This proves that $-(\Wc-\e)$ is accretive on $(\HH,\nr\cdot_{\HH,\nu})$.

Now for $\z \in \C^+$ we can check that the operator
\begin{equation*}
\Rc_\Wc(\z)
=
\begin{pmatrix}
-R(i\z) (aw + \z w) & - R(i\z)\\
w - w R(i\z) (\z aw + \z^2 w) & -\z w R(i\z)
\end{pmatrix}
\end{equation*}
is bounded on $\HH$, with $\Ran(\Rc_\Wc(\z)) \subset \Dom(\Wc)$. Moreover, $\Rc_{\Wc}(\z) (\Wc-\z) = \Id_{\Dom(\Wc)}$ and $(\Wc-\z) \Rc_{\Wc}(\z) = \Id_{\HH}$. This proves that $\z \in \rho(\Wc)$, and in particular, $-(\Wc-\e)$ is maximal accretive on $(\HH,\nr\cdot_{\HH,\nu})$. By the Lumer-Phillips Theorem, it generates a contractions semigroup on $(\HH,\nr\cdot_{\HH,\nu})$. And by equivalence of the norms, there exists $M > 0$ such that
\[
\forall t \geq 0, \quad \|e^{t(\Wc-\e)}\|_{\Lc(\HH)} \leq M.
\]
The conclusion follows.
\end{proof}

Proposition \ref{prop:Wc} ensures that for $F = (f,g) \in \HH$ the Cauchy problem
\begin{equation} \label{wave-Wc}
\begin{cases}
\partial_t U(t) = \Wc U(t), \quad t\geq 0,\\
U(0) = F,
\end{cases}
\end{equation}
has a solution $U \in C^0(\R_+;\HH)$. If moreover $F \in \Dom(\Wc)$, then $U \in C^0(\R_+;\Dom(\Wc)) \cap C^1(\R_+,\HH)$. Denoting by $u$ the first component of $U$ (the second being then $w \partial_t u$), we get in particular the following well-posedness result. We recall that by density it will be enough to work with $(f,g) \in \Sc \times \Sc$.

\begin{corollary} \label{cor:well-posed}
For $(f,g) \in \Sc\times \Sc$ the problem \eqref{wave} has a unique solution $u \in C^0(\R_+,H^2) \cap C^1(\R_+,H^1) \cap C^2(\R_+,L^2)$. Moreover, for $\e > 0$ there exists $C_\e > 0$ such that
\begin{equation*}
\forall t \geq 0, \quad \nr{u(t)}_{H^1}^2 + \nr{\partial_t u(t)}_{L^2}^2 \leq C_\e^2 e^{2t\e} \big( \nr{f}_{H^1}^2 + \nr{g}_{L^2}^2 \big).
\end{equation*}
\end{corollary}

We now turn to the proof of the estimates of Theorem \ref{th:loc-decay}. Let $(f,g) \in \Sc \times \Sc$ and let $u$ be the solution of \eqref{wave}. Let $\chi \in C^\infty(\R;[0,1])$ be equal to 0 on $]-\infty,1]$ and equal to 1 on $[2,+\infty[$.  For $\mu > 0$ and $t \in \R$ we set
\[
u_{\chi,\mu}(t) = \chi(t) e^{-t\mu} u(t).
\]
By Corollary \ref{cor:well-posed}, this defines a function in $\Sc(\R;L^2)$. By Fourier transform, for all $t \in \R$ we can write 
\[
u_{\chi,\mu}(t) = \frac 1 {2\pi} \int_\R e^{-it\tau} v_{\chi,\mu}(\tau) \diff \tau,
\]
where, for $\tau \in \R$,
\begin{equation} \label{def:v-chi}
v_{\chi,\mu}(\tau) = \int_\R e^{it\tau} u_{\chi,\mu}(t) \diff t, \quad \forall \tau \in \R.
\end{equation}

Given $\tau \in \R$, we write $z$ for $\tau + i\mu \in \C_+$. Then we have
\begin{align*}
    P(z) v_{\chi,\mu}(\tau) = F_{\chi,z}
\end{align*}
where
\begin{align*}
F_{\chi,z} & =\int_\R e^{itz} \big( aw \chi'(t) + 2 w \chi'(t) \partial_t + w \chi''(t) \big) u(t) \diff t \\&
=\int_\R e^{itz} \big( aw \chi'(t) - 2 iz w \chi'(t) - w \chi''(t) \big) u(t) \diff t.
\end{align*}
Finally, for $t \geq 2$,
\[
u(t) = e^{t\mu} u_{\chi,\mu}(t) = \frac 1 {2\pi} \int_\R e^{-itz} R(z) F_{\chi,z} \diff \tau.
\]

Now we separate  the contributions  of low  and high frequencies.  For this, we consider $\phi\in C^\infty(\R;[0,1])$ supported in $]-2,2[$ and equal to 1 on a neighborhood of $[-1,1]$. Then for $t \geq 2$, we set
\[
u_{\mathsf{low},\chi,\mu}(t) =  \frac 1 {2\pi} \int_\R \phi(\tau)e^{-itz} R(z) F_{\chi,z} \diff \tau,
\]
and
\[
u_{\mathsf{high},\chi,\mu}(t) =  \frac 1 {2\pi} \int_\R (1-\phi)(\tau)e^{-itz} R(z) F_{\chi,z} \diff \tau.
\]

For $\d \in \R$, we set $\HH^\d = H^{1,\d} \times L^{2,\d}$. We have the following result of propagation at finite speed.

\begin{proposition} \label{prop:prop-finite-speed}
Let $\d \geq 0$ and $T > 0$. There exists $C_{T,\d} > 0$ such that for $F = (f,g) \in \HH^\d$ and $t \in [0,T]$, we have
\[
\nr{u(t)}_{L^{2,\d}} \leq C_{T,\d} \nr{F}_{\HH^\d},
\]
where $u$ is the solution of \eqref{wave}.
\end{proposition}

\begin{proof}
\stepp  We set $\gamma = \max \big(C_G , C_w \big) \geq 1$, where $C_G$ and $C_w$ are the constants which appear in \eqref{eq:CGCw}. First assume that $F$ is supported in the ball $B(R)$ for some radius $R > 0$. We have
\begin{eqnarray*}
\lefteqn{\frac {\mathrm d}{\mathrm dt} \int_{\R^d \setminus B(R+\gamma t)} \big( \innp{G \nabla u(t)}{\nabla u(t)} + w \abs{\partial_t u(t)}^2 \big) \diff x}\\
&& = \int_{\R^d \setminus B(R+\gamma t)} 2 \Re \big( \innp{G \nabla u(t)}{\nabla \partial_t u(t)} + w \partial_t^2 u(t) \overline{\partial_t u(t)} \big) \diff x \\
&& \quad - \gamma  \int_{\abs x = R+\gamma t} \big( \innp{G\nabla u(t)}{\nabla u(t)} + w\abs{\partial_t u(t)}^2 \big) \diff \s(x).
\end{eqnarray*}
Since
\begin{eqnarray*}
\lefteqn{2 \Re\int_{\R^d \setminus B(R+\gamma t)}  \big( \innp{G \nabla u(t)}{\nabla \partial_t u(t)} + w \partial_t^2 u(t) \overline{\partial_t u(t)} \big) \diff x} \\
&& = - 2 \Re \int_{\R^d \setminus B(R+\gamma t)}  aw \abs{\partial_t u(t)}^2 \diff x - 2 \Re \int_{\abs x = R + \gamma t} \frac {x}{\abs {x}} \cdot (G\nabla u(t)) \overline{\partial_t u(t)} \diff \s(x) \\
&& \leq  \int_{\abs x = R + \gamma t} \big( \langle G G^{\frac 12}\nabla u(t),G^{\frac 12}\nabla u(t) \rangle  + \abs{\partial_t u(t)}^2 \big) \diff \s(x) \\
&& \leq \gamma \int_{\abs x = R+\gamma t} \big( \innp{G\nabla u(t)}{\nabla u(t)} + w\abs{\partial_t u(t)}^2 \big) \diff \s(x),
\end{eqnarray*}
we get
\[
\frac {\mathrm d}{\mathrm dt} \int_{\R^d \setminus B(R+\gamma t)} \big( \innp{G \nabla u(t)}{\nabla u(t)} + w \abs{\partial_t u(t)}^2 \big) \diff x \leq 0.
\]
Then
\[
\int_{\R^d \setminus B(R+\gamma t)} \big( \innp{G \nabla u(t)}{\nabla u(t)} + w \abs{\partial_t u(t)}^2 \big) \diff x = 0
\]
for all $t \geq 0$, since this is the case for $t = 0$ by assumption. We deduce that the solution $u(t)$ itself vanishes outside $B(R+\gamma t)$.

\stepp Now we consider a general $F \in \HH^\d$. We set $\Cc_0 = \set{x \in \R^d \, : \, \abs x \leq 2 \gamma T}$ and $\Cc_n = \set{x \in \R^d \, : \, 2^{n} \gamma T \leq \abs x \leq 2^{n+1} \gamma T}$ for $n \in \N^*$. We consider $\chi \in C_0^\infty(\R^d;[0,1])$ equal to 1 on $B(2 \gamma T)$ and supported in $B(3\gamma T)$. Then, for $n \in \N^*$, we consider $\chi_n : x \mapsto \chi(2^{-n} x)$, and we denote by $u_n$ the solution of \eqref{wave} with initial condition $(1-\chi_n)F$.
Let $n \geq 2$. By propagation at finite speed, since $(u-u_{n-2})(0)$ is supported in $B(3\cdot 2^{n-2} \g T)$, then $u(t)$ coincides with $u_{n-2}(t)$ on $\Cc_{n}$ for any $t \in [0, T]$.  Then, by Corollary \ref{cor:well-posed}, we have for $t \in [0, T]$ and $n \geq 2$
\begin{align*}
\nr{u(t)}_{L^2(\Cc_n)}^2
 = \nr{u_{n-2}(t)}_{L^2(\Cc_n)}^2
 \lesssim \nr{(1-\chi_{n-2}) F}_{\HH}^2
 \lesssim \sum_{k \geq n-1} \nr{F}_{H^1(\Cc_k) \times L^2(\Cc_k)}^2.
\end{align*}
Moreover, $\nr{u(t)}_{L^2(\Cc_0\cup \Cc_1)} \lesssim \nr{u(t)} \lesssim \nr{F}_\HH$, so
\begin{align*}
\nr{u(t)}_{L^{2,\d}}^2
& \lesssim \sum_{n \in \N} 2^{2(n+1) \d} \nr{u(t)}_{L^2(\Cc_n)}^2\\
& \lesssim \nr{F}_{\HH}^2 + \sum_{n \geq 2} 2^{2(n+1) \d} \sum_{k \geq n-1}  \nr{F}_{H^{1}(\Cc_k) \times L^{2}(\Cc_k)}^2\\
& \lesssim \nr{F}_{\HH}^2 + \sum_{n \geq 2} 2^{2(n+1) \d} \sum_{k \geq n-1} 2^{-2k\d} \nr{F}_{H^{1,\d}(\Cc_k) \times L^{2,\d}(\Cc_k)}^2\\
& \lesssim \nr{F}_{\HH}^2 + \sum_{k \in \N^*} 2^{-2k\d} \nr{F}_{H^{1,\d}(\Cc_k) \times L^{2,\d}(\Cc_k)}^2 \sum_{n \leq k+1} 2^{2(n+1) \d} \\
& \lesssim \nr{F}_{\HH^\d}^2 .
\end{align*}
The proof is complete.
\end{proof}

For the contribution of high frequencies we have the following estimate.

\begin{proposition}\label{prop:time-dec-high-frequ}
Let $n \in \N$ and $\d > n + \frac 12$. There exists $C > 0$ independent of $f,g \in \Sc$ such that for all $\mu > 0$ and $t \geq 1$, we have
\[
\nr{u_{\mathsf{high},\chi,\mu}(t)}_{L^{2,-\d}} \leq C e^{t\mu} t^{-n} \big( \nr{f}_{H^{1,\d}} + \nr{g}_{L^{2,\d}} \big).
\]
\end{proposition}

\begin{proof}
As above, we write $z$ for $\tau + i \mu$, $\tau \in \R$. By integrations by parts, we have
\begin{align*}
(it)^n u_{\mathsf{high},\chi,\mu}(t)
= \frac 1 {2\pi} \int_\R e^{-itz}  \partial_\tau^n \big((1-\phi)(\tau) R(z) F_{\chi,z} \big) \diff \tau.
\end{align*}
This can be written as a sum of terms of the form
\[
T_{\mu,n_1,n_2,n_3}(t) = \frac 1 {2\pi} \int_\R e^{-itz}   (1-\phi)^{(n_1)}(\tau) R^{(n_2)}(z) \partial_z^{n_3}F_{\chi,z}  \diff \tau,
\]
where $n_1 + n_2 + n_3 = n$. For such a term we have, by Theorem \ref{th:high-freq},
\begin{align*}
\nr{T_{\mu,n_1,n_2,n_3}(t)}_{L^{2,-\d}}
 \lesssim e^{t\mu} \int_{\abs \tau \geq 1} \abs\tau\inv \nr{\partial_z^{n_3}F_{\chi,z}}_{L^{2,\d}} \diff \tau 
 \lesssim e^{t\mu} \left(\int_{\R}  \nr{\partial_z^{n_3}F_{\chi,z}}_{L^{2,\d}}^2 \diff \tau \right)^{\frac 12}.
\end{align*}
By the Plancherel Theorem and Proposition \ref{prop:prop-finite-speed},
\[
\int_{\R} \nr{\partial_z^{n_3}F_{\chi,z}}_{L^{2,\d}}^2 \diff \tau \lesssim \int_{1 \leq t \leq 2} \nr{u(t)}_{L^{2,\d}}^2 \diff t \lesssim \nr{f}_{H^{1,\d}} ^2 + \nr{g}_{L^{2,\d}}^2.
\]
The conclusion follows.
\end{proof}

It was convenient for the contribution of high frequencies to introduce the cut-off function $\chi(t)$. Then the initial condition was somehow replaced by a quantity $F_{\chi,z}$ which depends on the values of the solution for $t\leq 2$. However, to get a precise asymptotic profile, it is now simpler to really work with $f$ and $g$. Then we set
\[
F_z = (aw-izw) f + wg.
\]
Notice that we formally recover $F_z$ from $F_{\chi,z}$ with $\chi$ replaced by $\1_{\R_+}$. In particular, we have
\begin{equation} \label{def:v}
P(z) v_\mu(\tau) = F_z, \quad \text{where} \quad  v_\mu(\tau) = \int_0^{+\infty} e^{it\tau} e^{-t\mu} u(t) \diff t.
\end{equation}
Then we set
\[
u_{\mathsf{low},\mu}(t) =  \frac 1 {2\pi} \int_\R  \phi(\tau)e^{-itz} R(z) F_{z} \diff \tau.
\]

We check that the difference between $u_{\mathsf{low},\chi,\mu}(t)$ and $u_{\mathsf{low},\mu}(t)$ is irrelevant for our purpose.

\begin{proposition} \label{prop:comparaison-ulow-1}
Let $m \in \N$. There exists $C > 0$ independent of $f,g \in \Sc$ such that for all $\mu > 0$ and $t\geq 2$, we have
\[
\nr{u_{\mathsf{low},\chi,\mu}(t)- u_{\mathsf{low},\mu}(t)}_{L^2} \leq C  e^{t\mu}\pppg t^{-m} \big(\nr{f}_{H^1} + \nr{g}_{L^2} \big).
\]
\end{proposition}

\begin{proof}
By integrations by parts, we have
\[
(it)^m \big(u_{\mathsf{low},\chi,\mu}(t)- u_{\mathsf{low},\mu}(t)\big) = \frac{1}{2\pi} \int_\R e^{-itz} \partial_\tau^m \big(\phi(\tau)\big(v_{\chi,\mu}(\tau) - v_{\mu}(\tau) \big) \big) \diff \tau.
\]
Let $\nu \in \Ii 0 m$. By \eqref{def:v-chi} and \eqref{def:v}, 
\[
\partial_\tau^\nu \big(v_{\chi,\mu}(\tau) - v_{\mu}(\tau) \big) = \int_0^2 (is)^\nu e^{is\tau} e^{-s\mu} (\chi(s)-1) u(s) \diff s.
\]
Consequently, with Corollary \ref{cor:well-posed}
\[
\nr{\partial_\tau^\nu \big(v_{\chi,\mu}(\tau) - v_{\mu}(\tau) \big)}_{L^2}   \lesssim \sup_{0\leq s \leq 2} \nr{u(s)}_{L^2} \lesssim \nr f_{H^1}  + \nr g_{L^2},
\]
and the conclusion follows.
\end{proof}

Finally, we have proved that for $n \in \N$ and $\d > n + \frac 12$ there exists $C > 0$ such that for all $\mu > 0$ and $t \geq 0$
\[
\nr{u(t) - u_{\mathsf{low},\mu}(t)}_{L^{2,-\d}} \leq C e^{t\mu}  \pppg t^{-n} \big( \nr{f}_{H^{1,\d}} + \nr{g}_{L^{2,\d}} \big).
\]
With a possibly different constant $C$, we also have
\[
\nr{u(t) - u_{\mathsf{low},\mu}(t)}_{L^{2,-\d}} \leq C e^{t\mu}  \pppg t^{-n} \big( \nr{f}_{H^{1,\d}} + \nr{af+g}_{L^{2,\d}} \big).
\]
Similarly, given $(f_0,g_0) \in \HH^\d$ the solution $u_0(t)$ of the free wave equation \eqref{wave-free} is close to
\[
u_{0,\mathsf{low},\mu}(t)= \frac 1 {2\pi} \int_\R  \phi(\tau)e^{-itz} R_0(z) F_{0,z} \diff \tau,  \]
where  $F_{0,z}= g_0-izf_0$, in the sense that
\[
\nr{u_0(t) - u_{0,\mathsf{low},\mu}(t)}_{L^{2,-\d}} \leq C e^{t\mu} \pppg t^{-n} \big( \nr{f_0}_{L^{2,\d}} + \nr{af_0 + g_0}_{L^{2,\d}} \big).
\]
Now we can use Theorem \ref{th:low-freq} to compare $u_{\mathsf{low},\mu}(t)$ with $u_{0,\mathsf{low},\mu}(t)$. Notice that in Theorem \ref{th:loc-decay} we have chosen $f_0=fw$ and $g_0= wg+awf$ to have $F_z = F_{0,z}$. We use a classical argument to convert the resolvent estimates into time decay.

\begin{proposition} \label{prop:comparaison-ulow-2}
Let $\rho_1 \in ]0,\rho_0[$. Let $\d > d + \frac 52$. There exists $C > 0$ independent of $f,g \in \Sc$ such that for all $\mu > 0$ and $t\geq 2$, we have
\[
\nr{u_{\mathsf{low},\mu}(t)- u_{0,\mathsf{low},\mu}(t)}_{L^{2,-\d}} \leq C e^{t\mu} \left(\pppg t^{-d-\rho_1} \nr{f}_{L^{2,\d}} +  \pppg t^{1-d-\rho_1} \nr{af + g}_{L^{2,\d}} \right),
\]
where $(f_0,g_0) = (fw,wg+awf)$.
\end{proposition}

\begin{proof}
We set $h = af + g$. As above, we write $z$ for $\tau + i\mu$. Since $F_z = F_{0,z}$ we have for $t \geq 2$,
\[
u_{\mathsf{low},\mu}(t)-u_{0,\mathsf{low},\mu}(t) =  \frac 1 {2\pi} \int_\R \phi(\tau)e^{-itz} \big( R(z) -R_0(z)\big) F_z \diff \tau.
\]
By integrations by parts, we have
\begin{equation}
\label{R1}
(it)^{d-1} \big( u_{\mathsf{low},\mu}(t)-u_{0,\mathsf{low},\mu}(t)\big) =  \frac 1 {2\pi} \int_\R e^{-itz} \partial_\tau^{d-1} \big(\theta^1_{\mu}(z) + \theta^2_\mu(z) \big) \diff \tau,
\end{equation}
where
\[
\theta^1_{\mu}(z)=\phi(\tau) \big(R(z) -R_0(z)\big) wh, \quad \theta^2_{\mu}(z)= -i\phi(\tau) \big(R(z) -R_0(z)\big) zwf.
\]
By Theorem \ref{th:low-freq} we have, for $\nu \in \{d-1,d\}$,
\[
\nr{\partial_\tau^\nu\theta^1_{\mu}(z)}_{L^{2,-\d}}  \lesssim |\tau|^{d+\rho_1-\nu-2} \nr{h}_{L^{2,\d}}.
\]
Then, on the one hand,
\begin{align*}
\nr{\int_{-t^{-1}}^{t^{-1}} e^{-itz} \partial_\tau^{d-1}\theta^1_{\mu}(z) \diff \tau}_{L^{2,-\d}}
\lesssim  e^{t\mu} \nr{h}_{L^{2,\d}} \int_{-t^{-1}}^{t^{-1}} |\tau|^{\rho_1-1} \diff \tau
\lesssim  e^{t\mu}  t^{-\rho_1} \nr{h}_{L^{2,\d}}.
\end{align*}
On the other hand, with another integration by parts, we get
\begin{eqnarray*}
\lefteqn{t \nr{\int_{|\tau|\geq t^{-1}} e^{-itz}  \partial_\tau^{d-1}\theta^1_{\mu}(z) \diff \tau}_{L^{2,-\d}}} \\
&& \leq e^{t\mu}  \nr{\partial_\tau^{d-1}\theta^1_{\mu}(-t^{-1} + i\mu)}_{L^{2,\d}} + e^{t\mu} \nr{\partial_\tau^{d-1}\theta^1_{\mu}(t^{-1}+i\mu)
}_{L^{2,-\d}}\\
&& \quad + e^{t\mu}  \int_{t^{-1}\leq |\tau| \leq 2}  \nr{\partial_\tau^{d}\theta^1_{\mu}(z)}_{L^{2,-\d}}  \diff \tau  \\
&&  \lesssim  e^{t\mu} t^{1-\rho_1}  \nr{h}_{L^{2,\d}}
+e^{t\mu} \nr{h}_{L^{2,\d}} \int_{t^{-1}\leq |\tau| \leq 2}
|\tau|^{\rho_1-2} \diff \tau \\
&&\lesssim  e^{t\mu}     t^{1-\rho_1} \nr{h}_{L^{2,\d}}.
\end{eqnarray*}
Finally,
\begin{equation} \label{R2}
\nr{\int_\R e^{-itz} \partial_\tau^{d-1}\theta^1_{\mu}(z) \diff \tau}_{L^{2,-\d}}
\lesssim  e^{t\mu}  t^{-\rho_1} \nr{h}_{L^{2,\d}}.
\end{equation}
For $\th^2_\mu$ we have an additional power of $z$, so we do another integration by parts to write
\begin{equation} \label{R3a}
it \int_\R e^{-itz} \partial_\tau^{d-1}  \theta^2_\mu(z) \diff \tau = \int_\R e^{-itz} \partial_\tau^{d}  \theta^2_\mu(z) \diff \tau.
\end{equation}
Then we can proceed as above. Since for $\nu \in \{d,d+1\}$, we have
\[
\nr{\partial_\tau^\nu\theta^2_{\mu}(z)}_{L^{2,-\d}}  \lesssim |\tau|^{d+\rho_1-\nu-1} \nr{f}_{L^{2,\d}},
\]
we get
\begin{equation} \label{R3}
\nr{\int_\R e^{-itz} \partial_\tau^d\theta^2_{\mu}(z) \diff \tau}_{L^{2,-\d}}
\lesssim  e^{t\mu}  t^{-\rho_1} \nr{f}_{L^{2,\d}}.
\end{equation}
The conclusion follows from \eqref{R1}, \eqref{R2}, \eqref{R3a} and \eqref{R3}.
\end{proof}

Now we can conclude the proof of Theorem \ref{th:loc-decay}.

\begin{proof}[Proof of Theorem \ref{th:loc-decay}]
With Propositions \ref{prop:time-dec-high-frequ}, \ref{prop:comparaison-ulow-1} and \ref{prop:comparaison-ulow-2}, we have proved that there exists $C > 0$ independent of $\mu > 0$, $f \in \Sc$, $g \in \Sc$ and $t \geq 0$ such that
\[
\nr{u(t) - u_0(t)}_{L^{2,-\d}} \leq C e^{t\mu} \big(\pppg t^{-d-\rho_1} \nr{f}_{H^{1,\d}} +   \pppg t^{1-d-\rho_1} \nr{af + g}_{L^{2,\d}}\big).
\]
By density, the same inequality holds for any $f \in H^{1,\d}$ and $g \in L^{2,\d}$. Finally, since the left-hand side does not depend on $\mu$, it only remains to let $\mu$ go to 0 to conclude.
\end{proof}

\section{Strategy of the proof for low frequencies} \label{low-freq}

In this section, we give the main arguments for the proof of Theorem \ref{th:low-freq}. Some intermediate results will be proved in the following sections. We set
\[
\DD = \set{\z \in \C \, : \, \abs \z \leq 1}, \quad \DD_+ = \DD \cap \C_+,
\]
so that the estimates of Theorem \ref{th:low-freq} concern $z \in \DD_+$.

We have to compare the resolvents $R(z)$ and $R_0(z)$ (see \eqref{def:R}, \eqref{def:R0} and Proposition \ref{def:Pz}). For this, we use the resolvent identity
\begin{equation}
R(z) - R_0(z) = R(z) \th(z) R_0(z),
\end{equation}
where
\begin{equation}
\th(z) =  (\D_G - \D) +izaw +z^2(w-1).
\end{equation}

We actually compare the difference of the derivatives of these resolvents. For $n \in \N$, the difference $R^{(n)}(z)- R_0^{(n)}(z)$ is a sum of terms of the form
\[
R^{(n_1)}(z) \th^{(n_2)}(z) R^{(n_3)}_0(z), \quad n_1 + n_2 + n_3 = n.
\]
Notice that $\th^{(n_2)}(z) = 0$ for $n_2 \geq 3$. For $\sigma \in \{0,1,2\}$, we set $\th_\s(z) = \th^{(2-\sigma) }(z)$, so that
\begin{equation} \label{def:theta-sigma}
\th_0(z) = 2(w-1), \quad \th_1(z)= iaw + 2z(w-1), \quad \th_2(z) = \th(z).
\end{equation}
The motivation for this notation is that in suitable spaces $\th_\s(z)$ will be of size $\abs z^{\sigma + \rho}$ for $\rho \in ]0,\rho_0[$ (see Proposition \ref{prop:Pa-Po-comm} below for a precise statement), so this parameter $\sigma$ will appear in all the estimates below (in particular in Propositions \ref{prop:low-freq} to \ref{prop:regularity-Rz-A}).

We have 
\[
R'(z) = R(z) (iaw + 2zw) R(z).
\]
Setting
\begin{equation} \label{def:gamma01}
\g_1(z) = iaw+2zw, \quad \g_0(z) = \g_1'(z) = 2w,
\end{equation}
we can check by induction on $n_1 \in \N$ that the derivative $R^{(n_1)}(z)$ can be written as a sum of terms for the form
\begin{equation} \label{eq:der-Rz}
R(z)\g_{j_1}(z)R(z)...\g_{j_k}(z)R(z)
\end{equation}
where $j_1,\dots,j_k \in \lbrace 0,1\rbrace$ are such that
\begin{equation} \label{eq:der-k-n}
2k - (j_1+\dots+j_k)=n_1.
\end{equation}
The same applies to $R_0^{(n_3)}(z)$, with $R(z)$, $\g_1(z)$ and $\g_0(z)$ replaced by $R_0(z)$, $\g_1^0(z) = 2z$ and $\g_0^0(z) =  2$.

We will introduce $\g_2(z)$ and $\g_2^0(z)$ in \eqref{def:gamma2} below.
For $k \in \N$, $\jj = (j_1,\dots,j_k) \in \{0,1,2\}^k$ and $z,z' \in \DD_+$, we set
\[
\Rc_{k,\jj}(z',z)=R(z')\gamma_{j_1}(z) R(z')... \gamma_{j_k}(z) R(z'),
\]
\[
\Rc^0_{k,\jj}(z',z)=R_0(z')\gamma^0_{j_1}(z) R_0(z')... \gamma^0_{j_k}(z) R_0(z')
\]
(with the natural convention that this is just $R(z')$ or $R_0(z')$ for $k = 0$), and
\[
m(k,\jj) = 2(k+1) - \abs \jj = 2(k+1) - (j_1+\dots+j_k)  \geq 2.
\]
%
Notice for further reference that when $\jj \in \{0,1\}^k$, we have
\begin{equation} \label{eq:m-geq-k+2}
m(k,\jj) \geq k  +2.
\end{equation}
We write $\Rc_{k,\jj}(z)$ for $\Rc_{k,\jj}(z,z)$. In particular, each term \eqref{eq:der-Rz} is of the form $\Rc_{k,\jj}(z)$ with $k \in \N$, $\jj \in \{0,1\}^k$ and
\[
m(k,\jj) = n_1+2.
\]
The reason for introducing $m(k,\jj)$ is that for suitable $z,z'$ ($\abs z = \abs {z'} \sim \Im(z'))$ then $\Rc(z',z)$ will be of size $\abs z^{-m(k,\jj)}$ (see Proposition \ref{prop:reg-T} below). Then this quantity will also appear in all the estimates below (see again Propositions \ref{prop:low-freq} to \ref{prop:regularity-Rz-A}).\\

Finally, for $z \in \DD_+$ the difference $R^{(n)}(z) - R_0^{(n)}(z)$ can be written as a sum of operators of the form
\[
\Rc_{k,\jj}(z) \th_\sigma(z) \Rc^0_{k_0,\jj_0}(z),
\]
with $\sigma \in \{0,1,2\}$, $k,k_0 \in \N$, $\jj \in \{0,1\}^k$ and $\jj_0 \in \{0,1\}^{k_0}$ such that
\[
m(k,\jj) - \s + m(k_0,\jj_0) = n + 2.
\]

Then Theorem \ref{th:low-freq} is a consequence of the following estimate.

\begin{proposition} \label{prop:low-freq}
Let $\rho_1 \in ]0,\rho_0[$. Let $\sigma \in \{0,1,2\}$, $k,k_0 \in \N$, $\jj \in \{0,1\}$ and $\jj_0 \in \{0,1\}^{k_0}$. Let $\d > m(k,\jj) + m(k_0,\jj_0) - \sigma - \frac 12$. Then there exists $C > 0$ such that for $z \in \DD_+$, we have
\begin{equation} \label{estim-RR_0}
\nr{\pppg x^{-\d} \Rc_{k,\jj}(z) \th_\sigma(z) \Rc^0_{k_0,\jj_0}(z) \pppg x^{-\d}}_{\Lc(L^2)} \leq C \abs z^{\min(d - m(k,\jj)  - m(k_0,\jj_0) + \sigma + \rho_1,0)}.
\end{equation}
\end{proposition}

The main ingredient to prove such resolvent estimates near the real axis is the Mourre commutators method. In particular, for the resolvent of a Schr\"odinger-type operator, one usually uses the selfadjoint generator of dilations, defined by
\begin{equation} \label{def:A}
A = - \frac {x \cdot i\nabla + i\nabla \cdot x} 2 = - \frac {id}2 - x \cdot i\nabla.
\end{equation}
Its domain is the set of $u \in L^2$ such that ${(x \cdot \nabla)u} \in L^2$ in the sense of distributions.

In Section \ref{sec:Mourre}, we will deduce from the abstract commutators method the following result.

\begin{proposition}
\begin{enumerate}[\rm (i)]
\item Let $k \in \N$, $\jj \in \{0,1\}^k$ and $\d > k + \frac 12$. There exists $C > 0$ such that for $z \in \DD_+$, we have
\begin{equation} \label{estim-Mourre-wave-1}
\nr{\pppg {\A}^{-\d} \Rc_{k,\jj}(z) \pppg {\A}^{-\d}}_{\Lc(L^2)} \leq C \abs z^{-m(k,\jj)}.
\end{equation}
%
\item Let $\rho \in [0,\rho_0[$. Let $k,k_0 \in \N$, $\jj \in \{0,1\}^k$, $\jj_0 \in \{0,1\}^{k_0}$ and $\d > k + k_0 + \frac 3 2$. Let $\sigma \in \{0,1,2\}$. There exists $C > 0$ such that, for $z \in \DD_+$, 
\begin{equation} \label{estim-Mourre-wave-2}
\nr{\pppg {\A}^{-\d} \Rc_{k,\jj}(z) \th_\sigma(z) \Rc^0_{k_0,\jj_0}(z) \pppg {\A}^{-\d}}_{\Lc(L^2)} \leq  C \abs z^{-m(k,\jj) + (\sigma+\rho)  - m(k_0,\jj_0)}.
\end{equation}
\end{enumerate}\label{prop:Mourre}
\end{proposition}

These estimates do not directly give Proposition \ref{prop:low-freq}. The second ingredient is the fact that one can use the elliptic regularity given by the resolvents and the decay of the weights $\pppg x^{-\d}$ to get some smallness for low frequencies (this will be detailed in Section \ref{sec:hardy-ineq}). This is also the reason why the decay of the coefficients in $\th(z)$ gives extra smallness for the difference $R(z)-R_0(z)$ compared to the estimates for each resolvent alone.\\

To get elliptic regularity, we replace the resolvent $R(z)$ by $R(i\abs z)$. For $z \in \DD_+$ and $r = \abs z$, we have the resolvent identities
\begin{equation} \label{eq:res-identity}
R(z) - R(ir)=R(ir) \gamma_2(z) R(z)= R(z)\gamma_2(z) R(ir)
\end{equation}
and
\begin{equation} \label{eq:res0-identity}
R_0(z) - R_0(ir)=R_0(ir)\gamma_2^0(z)R_0(z)= R_0(z)\gamma_2^0(z)R_0(ir),
\end{equation}
where
\begin{equation} \label{def:gamma2}
\gamma_2(z)=(r+iz)aw+(z^2 + r^2)w  \quad \text{and} \quad  \gamma_2^0(z)=z^2+r^2
\end{equation}
 (notice that the factors commute in \eqref{eq:res0-identity} but not in \eqref{eq:res-identity}).

In order to get elliptic regularity, the idea is to replace each term in \eqref{eq:der-Rz} by a sum of terms with resolvents $R(ir)$ only, or with as many factors $R(ir)$ as needed. The proof of the following lemma is given in Appendix \ref{sec:additional-proofs}.

\begin{lemma} \label{lem:T-N-resolventes}
Let $k \in \N$, $\jj \in \{0,1\}^k$, $z \in \DD_+$ and $r = \abs z$. Let $N \in \N$. Then $\Rc_{k,\jj}(z)$ can be written as a sum of terms of the form
\begin{equation} \label{eq:T-N-res-1}
\Rc_{\k,\ll}(ir,z)
\end{equation}
for some $\k \in \{k,\dots,N\}$ and $\ll \in \{0,1,2\}^\k$ such that
\begin{equation} \label{eq:egalite-m}
m(\k,\ll) = m(k,\jj),
\end{equation}
or of the form
\begin{equation} \label{eq:T-N-res-2}
\Rc_{N,\jj_1}(ir,z) \gamma_{\ell}(z) \Rc_{k_2,\jj_2}(z),
\end{equation}
with $k_2 \in \{0,\dots ,k\}$, $\jj_1 \in \{0,1,2\}^{N}$, $\jj_2 \in \{0,1\}^{k_2}$ and $\ell \in \{0,1,2\}$ such that
\begin{equation} \label{eq:egalite-m-m}
m(N,\jj_1) + m(k_2,\jj_2) - \ell = m(k,\jj).
\end{equation}
\end{lemma}

Similarly, $\Rc^0_{k_0,\jj_0}(z)$ can be written as a sum of terms of the form
\begin{equation} \label{eq:T-N-res0-1}
\Rc^0_{\k_0,\ll_0}(ir,z)
\end{equation}
for some $\k_0 \in \{k_0,\dots,N\}$ and $\ll_0 \in \{0,1,2\}^{\k_0}$ such that
\begin{equation} \label{eq:egalite-m0}
m(\k_0,\ll_0) = m(k_0,\jj_0),
\end{equation}
or of the form
\begin{equation} \label{eq:T-N-res0-2}
\Rc^0_{k_1,\ll_1}(z) \gamma_{\ell_0}^0(z) \Rc^0_{N,\ll_2}(ir,z),
\end{equation}
with $k_1 \in \{0,\dots ,k_0\}$, $\ll_1 \in \{0,1\}^{k_1}$, $\ll_2 \in \{0,1,2\}^{N}$ and $\ell_0 \in \{0,1,2\}$ such that
\begin{equation} \label{eq:egalite-m0-m0}
m(k_1,\ll_1) + m(N,\ll_2) - \ell_0 = m(k_0,\jj_0).
\end{equation}

Thus, we will prove \eqref{estim-RR_0} with $\Rc_{k,\jj}(z)$ replaced by terms of the form \eqref{eq:T-N-res-1} or \eqref{eq:T-N-res-2}, and $\Rc^0_{k_0,\jj_0}(z)$ replaced by terms of the form and \eqref{eq:T-N-res0-1} or \eqref{eq:T-N-res0-2}.

For the contributions \eqref{eq:T-N-res-1} and \eqref{eq:T-N-res0-1}, which only involve resolvents of the form $R(ir)$ or $R_0(ir)$, we use the following estimate.

\begin{proposition} \label{prop:regularity-Rz}
Let $\rho \in [0,\rho_0[$. Let $k,k_0 \in \N$ , $\jj \in \{0,1,2\}^k$ and $\jj_0 \in \{0,1,2\}^{k_0}$. Let $s_1,s_2 \in \big[0,\frac d 2\big[$, $\d_1 > s_1$ and $\d_2 > s_2$.
There exists $C > 0$ such that for $z \in \DD_+$ and $r = \abs z$, we have 
\begin{equation*}
\nr{\pppg x^{-\d_1} \Rc_{k,\jj}(ir,z) \pppg x^{-\d_2}}_{\Lc(L^2)}
\leq C \abs z^{\min(s_1 + s_2 - m(k,\jj),0)},
\end{equation*}
and
\begin{equation*}
\nr{\pppg x^{-\d_1} \Rc_{k,\jj}(ir,z)\th_\sigma(z) \Rc^0_{k_0,\jj_0}(ir,z) \pppg x^{-\d_2}}_{\Lc(L^2)}
\leq C \abs z^{\min(s_1 + s_2 - m(k,\jj) + (\sigma+ \rho) - m(k_0,\jj_0),0)}.
\end{equation*}
\end{proposition}

For the terms \eqref{eq:T-N-res-2} and \eqref{eq:T-N-res0-2} which still involve resolvents of the form $R(z)$ or $R_0(z)$, we will apply the commutators method (see Proposition \ref{prop:Mourre}). Then we use the elliptic regularity to get some smallness and to compensate for the derivatives which appear in the conjugate operator $\A$. More precisely, we need the following estimates.

\begin{proposition} \label{prop:regularity-Rz-A}
Let $\rho \in [0,\rho_0[$. Let $s \in \big[0,\frac d 2 \big[$ and $\d > s$. Let $\s \in \{0,1,2\}$. Let $k,N \in \N$,  $\jj , \jj_0 \in \{0,1,2\}^N$ and $\jj_1,\jj_2 \in \{0,1,2\}^k$. Let $\ell,\ell_0 \in \{0,1,2\}$.
There exist $N_0 \in \N$ and $C > 0$ such that if $N \geq N_0$ then for $z \in \DD_+$ and $r = \abs z$, we have
\begin{equation}
    \label{eq:xRA-Schro-1}
\big\| \pppg x^{-\d}  \Rc_{N, \jj}(ir,z)\gamma_{\ell}(z) \pppg {\A}^\d \big\|_{\Lc(L^2)} \leq  C \abs z^{s-m(N,\jj)+\ell},
\end{equation}
\begin{equation}
\label{eq:xRA-Schro-2}
\big\| \pppg x^{-\d} \Rc_{k,\jj_1}(ir,z)\th_\sigma(z) \Rc^0_{N,\jj_0}(ir,z)\gamma^0_{\ell_0}(z)\pppg {\A}^\d\big\|_{\Lc(L^2)}  \leq  C \abs z^{s  -m(k,\jj_1)+(\s+ \rho)-m(N,\jj_0)+ \ell_0},
\end{equation}
\begin{equation}
\label{eq:xRA-Schro-3}
\big\|\pppg {\A}^\d \gamma^0_{\ell_0}(z) \Rc^0_{N,\jj_0}(ir,z) \pppg x^{-\d}\big\|_{\Lc(L^2)}  \leq  C \abs z^{s-m(N,\jj_0)+\ell_0},
\end{equation}
\begin{equation}
\label{eq:xRA-Schro-4}
\big\|\pppg {\A}^\d  \gamma_{\ell}(z) \Rc_{N,\jj}(ir,z)\th_\sigma(z)\Rc^0_{k,\jj_2}(ir,z) \pppg x^{-\d}\big\|_{\Lc(L^2)}  \leq  C \abs z^{s-m(N,\jj)+(\s+\rho)-m(k,\jj_2)+\ell}.
\end{equation}
\end{proposition}

We will prove propositions \ref{prop:regularity-Rz} and \ref{prop:regularity-Rz-A} in Section \ref{sec:elliptic-regularity}. Now we have all the ingredients to prove Proposition \ref{prop:low-freq} (and hence Theorem \ref{th:low-freq}).

\begin{proof}[Proof of Proposition \ref{prop:low-freq}, assuming Propositions \ref{prop:Mourre}, \ref{prop:regularity-Rz} and \ref{prop:regularity-Rz-A}]
Let $\rho \in ]\rho_1,\rho_0[$. Let \( z \in \DD_+ \) and $r = \abs z$. As said above, we apply Lemma \ref{lem:T-N-resolventes} to replace $\Rc_{k,\jj}(z)$ and $\Rc^0_{k_0,\jj_0}(z)$ by terms of the form \eqref{eq:T-N-res-1} or \eqref{eq:T-N-res-2} and \eqref{eq:T-N-res0-1} or \eqref{eq:T-N-res0-2}, respectively. This gives four different forms of contributions to estimate.

We begin with the case where $\Rc_{k,\jj}(z)$ and $\Rc^0_{k_0,\jj_0}(z)$ are replaced by terms of the form \eqref{eq:T-N-res-1} and \eqref{eq:T-N-res0-1}, respectively. We apply Proposition \ref{prop:regularity-Rz} with
\begin{equation} \label{eq:choix-s}
s_1 = s_2 = \dfrac{1}{2} \min \big( d + \rho_1 - \rho, m(k,\jj) +m(k_0,\jj_0) - \s \big) < \frac d 2,
\end{equation}
and $\d_1 = \d_2 = \d$. By assumption we have $\d > m(k,\jj) +m(k_0,\jj_0) - \s - \frac 12$, so in particular
\[
\d > \frac 12 \big( m(k,\jj) +m(k_0,\jj_0)-\s \big).
\]
With \eqref{eq:egalite-m} and \eqref{eq:egalite-m0}, Proposition \ref{prop:regularity-Rz} gives
\begin{equation} \label{eq: R th R_O}
\begin{aligned}
\nr{\pppg x^{-\d} \Rc_{\k,\ll}(ir,z)\th_\sigma(z) \Rc^0_{\k_0,\ll_0}(ir,z) \pppg x^{-\d}}_{\Lc(L^2)}
&\lesssim \abs z^{\min(s_1 + s_2 -m(\k,\ll) + (\sigma + \rho) - m(\k_0,\ll_0),0)}\\
&\lesssim \abs z^{\min(s_1 + s_2 -m(k,\jj) + (\sigma + \rho) - m(k_0,\jj_0),0)}\\
& \lesssim \abs z^{\min(d+ \rho_1-m(k,\jj)+ \sigma  - m(k_0,\jj_0),0)}.
\end{aligned}
\end{equation}
Thus we have proved \eqref{estim-RR_0} for a contribution with factors \eqref{eq:T-N-res-1} and \eqref{eq:T-N-res0-1}.

Then we consider the case where $\Rc_{k,\jj}(z)$ and $\Rc^0_{k_0,\jj_0}(z)$ are replaced by terms of the form \eqref{eq:T-N-res-2} and \eqref{eq:T-N-res0-2}, respectively. Notice that by \eqref{eq:m-geq-k+2}, we have
\[
\d > m(k,\jj) + m(k_0,\jj_0) - \s - \frac 12 \geq k + k_0 +\frac 32 \geq k_2 + k_1 + \frac 32,
\]
so \eqref{estim-Mourre-wave-2} applies to $\Rc_{k_2,\jj_2}(z)  \th_\sigma(z)  \Rc^0_{k_1,\ll_1}(z)$. We choose $s$ as in \eqref{eq:choix-s}. We can assume that $N$ is greater that $N_0$ given in Proposition \ref{prop:regularity-Rz-A}. We write
\begin{align*}
\Nc :=
& \nr{\pppg x^{-\d} \Rc_{N,\jj_1}(ir,z) \gamma_{\ell}(z) \Rc_{k_2,\jj_2}(z)
 \th_\sigma(z)  \Rc^0_{k_1,\ll_1}(z) \gamma_{\ell_0}^0(z) \Rc^0_{N,\ll_2}(ir,z)\pppg x^{-\d}}_{\Lc(L^2)}\\
& \leq
\nr{\pppg x^{-\d} \Rc_{N,\jj_1}(ir,z) \gamma_{\ell}(z) \pppg {\A}^{\d}}_{\Lc(L^2)}\\
& \times
\nr{\pppg{\A}^{-\d} \Rc_{k_2,\jj_2}(z)  \th_\sigma(z)  \Rc^0_{k_1,\ll_1}(z) \pppg{\A}^{-\d}}_{\Lc(L^2)}\\
& \times \nr{\pppg{\A}^{\d}\gamma_{\ell_0}^0(z) \Rc^0_{N,\ll_2}(z) \pppg x^{-\d}}_{\Lc(L^2)}.
\end{align*}
By \eqref{eq:xRA-Schro-1}, \eqref{estim-Mourre-wave-2} and \eqref{eq:xRA-Schro-3}, we obtain with \eqref{eq:egalite-m-m} and \eqref{eq:egalite-m0-m0}
\begin{align*}
\Nc
& \lesssim \abs z^{s-m(N,\jj_1)+ \ell} \abs z^{-m(k_2,\jj_2) + (\s + \rho) - m(k_1,\mathsf l_1)} \abs z^{s-m(N,\ll_2)+ \ell_0}\\
& \lesssim \abs z^{2s-m(k,\jj) + (\s +\rho) - m(k_0,\jj_0)}\\
& \lesssim \abs z^{\min(d + \rho_1 -m(k,\jj)+ \s  - m(k_0,\jj_0) ,0)}.
\end{align*}

Now we consider the case where $\Rc_{k,\jj}(z)$ and $\Rc^0_{k_0,\jj_0}(z)$ are replaced by terms of the form \eqref{eq:T-N-res-1} and \eqref{eq:T-N-res0-2} respectively. We have to estimate an operator of the form
\[
\pppg x^{-\d} \Rc_{\k,\ll}(ir,z) \th_\sigma(z)  \Rc^0_{k_1,\ll_1}(z) \gamma_{\ell_0}^0(z) \Rc^0_{N,\ll_2}(ir,z)\pppg x^{-\d}.
\]
We cannot estimate directly this operator. We need more factors $R_0(ir)$, now on the left of $\Rc^0_{k_1,\ll_1}(z)$. Thus, as we did above, we write $\Rc^0_{k_1,\ll_1}(z)$ as a sum of terms of the form
\begin{equation} \label{eq:T-N-res00-1}
\Rc^0_{\k_0,\ll_0}(ir,z)
\end{equation}
for some $\k_0 \in \{k_1,\dots,N\}$ and $\ll_0 \in \{0,1,2\}^{\k_0}$ such that $m(\k_0,\ll_0) = m(k_1,\ll_1)$, or of the form
\begin{equation} \label{eq:T-N-res00-2}
\Rc^0_{N,\ll'_1}(ir,z) \gamma_{\ell'}^0(z) \Rc^0_{k_2,\ll'_2}(z),
\end{equation}
with $k_2 \in \{0,\dots ,k_1\}$, $\ll'_1 \in \{0,1,2\}^{N}$, $\ll'_2 \in \{0,1,2\}^{k_2}$ and $\ell' \in \{0,1,2\}$ such that $m(N,\ll'_1) + m(k_2,\ll'_2) - \ell' = m(k_1,\ll_1)$. A term of the form
\[
\pppg x^{-\d} \Rc_{\k,\ll}(ir,z) \th_\sigma(z)  \Rc^0_{\k_0,\ll_0}(ir,z) \gamma_{\ell_0}^0(z) \Rc^0_{N_2,\jj_2}(ir,z)\pppg x^{-\d}
\]
is estimated with Proposition \ref{prop:regularity-Rz} as before. On the other hand, with \eqref{eq:xRA-Schro-2}, \eqref{estim-Mourre-wave-2} and \eqref{eq:xRA-Schro-3},
\begin{eqnarray*}
\lefteqn{\nr{\pppg x^{-\d} \Rc_{\k,\ll}(ir,z) \th_\sigma(z)  \Rc^0_{N,\ll'_1}(ir,z) \gamma_{\ell'}^0(z) \Rc^0_{k_2,\ll'_2}(z) \gamma_{\ell_0}^0(z) \Rc^0_{N,\ll_2}(ir,z)\pppg x^{-\d}}_{\Lc(L^2)}}\\
&& \leq \nr{\pppg x^{-\d} \Rc_{\k,\ll}(ir,z) \th_\sigma(z)  \Rc^0_{N,\ll'_1}(ir,z) \gamma_{\ell'}^0(z) \pppg{\A}^{\d}} \nr{ \pppg{\A}^{-\d} \Rc^0_{k_2,\ll'_2}(z)\pppg{\A}^{-\d}} \\
&& \times \nr{\pppg{\A}^{\d} \gamma_{\ell_0}^0(z) \Rc^0_{N,\ll_2}(ir,z)\pppg x^{-\d}}\\
&& \lesssim \abs z^{s-m(\k,\ll) + (\s + \rho) - m(N,\ll'_1) +\ell'} \abs z^{-m(k_2,\ll'_2)} \abs z^{s-m(N_2,\ll_2)+ \ell_0}\\
&& \lesssim \abs z^{2s-m(k,\jj) + (\s +\rho) - m(k_0,\jj_0)}\\
&& \lesssim \abs z^{\min(d + \rho_1-m(k,\jj)+\s - m(k_0,\jj_0) ,0)}.
\end{eqnarray*}
The case where $\Rc_{k,\jj}(z)$ and $\Rc^0_{k_0,\jj_0}(z)$ are replaced by terms of the form \eqref{eq:T-N-res-2} and \eqref{eq:T-N-res0-1}, respectively, is similar.
\end{proof}

\section{The generalized Hardy inequality}\label{sec:hardy-ineq}

In the next two sections, we prove Propositions \ref{prop:Mourre}, \ref{prop:regularity-Rz} and \ref{prop:regularity-Rz-A}. For this, we fix $\rho \in ]0,\rho_0[$ and $\bar \rho \in ]\rho,\rho_0[$.\\

We first introduce suitable norms on the usual Sobolev spaces. For $r \in ]0,1]$ and $s \in \R$ we denote by $H^s_r$ the usual space $H^s$, endowed with the norm defined by 
\begin{equation*}
\nr{u}_{H^s_r} = \nr{\pppg {D_r}^s u}_{L^2}, \quad \text{where} \quad D_r = \sqrt{-\D}/r.
\end{equation*}
For $z \in \DD_+$, we set $H_z^s=H_{\abs z}^s$. Notice that for $\a \in \N^d$, we have
\begin{equation} \label{eq:derivee-alpha}
\nr{\partial_x^\a}_{\Lc(H_r^s,H_r^{s-\abs \a})} = r^{\abs \a}.
\end{equation}
These norms are adapted to the analysis of slowly oscillating functions, for which a derivative defines a "small" operator (these norms are useful when $r = \abs z$ is small).\\

One of the key points of the analysis is that multiplication by a coefficient decaying at infinity behaves like a derivative: it costs some regularity but it is small for low frequencies. For this, we generalize the usual Hardy inequality.

Let $d_0$ be a fixed integer  greater than $\frac d 2$. For $\kappa \geq 0$, we denote by $\Sc^{-\kappa}$ the set of smooth functions $\vf$ on $\R^d$ such that
\begin{equation} \label{def-nr-nu-N}
\nr{\vf}_{\Sc^{-\kappa}} =  \sup_{\abs \a \leq d_0}  \sup_{x \in \R^d} \big| \pppg x ^{\kappa +\abs \a} \partial^\a  \vf(x) \big| < +\infty.
\end{equation}

It is now classical that we can use the following result to convert decay at infinity and regularity into smallness for low frequencies. The following result is \cite[Prop.7.2]{BoucletRoy14} or \cite[Prop.3.1]{Royer24}.

\begin{proposition} \label{prop:dec-sob}
Let $s \in \big] -\frac d 2, \frac d 2\big[$ and $\kappa \geq 0$ be such that $s -\kappa \in \big] -\frac d 2, \frac d 2\big[$. Let $\eta > 0$. There exists $C \geq 0$ such that for $\vf \in \Sc^{-\kappa-\eta}$, $u \in H^s$ and $r \in ]0,1]$, we have
\[
 \nr{\vf u}_{H_r^{s-\kappa}} \leq C r ^\kappa \nr \vf_{\Sc^{-\kappa-\eta}} \nr u _{H_r^s}.
\]
\end{proposition}
\begin{remark} \label{rem:dec-sob}
If $\vf \in \Sc^{-\eta}$ for some $\eta > 0$, then for any $s \in \big]-\frac d 2,\frac d 2 \big[$ the multiplication by $(1+\vf)$ defines a bounded operator on $H_r^s$ uniformly in $r \in ]0,1]$.
\end{remark}

In particular, together with the elliptic regularity that will be given by the resolvents, the weights which appear on both sides in \eqref{estim-RR_0} will give positive powers of the frequency $\abs z$:

\begin{lemma} \label{lem:weight}
Let $s \in \big[0,\frac d 2 \big[$ and $\d > s$. There exists $C > 0$ such that for $r \in ]0,1]$, we have
\[
\|\pppg x^{-\d} \|_{\Lc(H_r^s,L^2)} \leq C \, r^s 
\quad \text{and} \quad
\|\pppg x^{-\d} \|_{\Lc(L^2,H_r^{-s})} \leq C \, r^s.
\]
\end{lemma}

We also use this generalized Hardy inequality to see that with the assumptions \eqref{hyp-swa} and \eqref{hyp-a} at infinity, the operator $P(z)$ (see \eqref{def:Pz}) is a perturbation of $P_0(z)$ is the following sense.

\begin{proposition} \label{prop:Pa-Po} 
Let $s \in \big]-\frac d 2, \frac d 2 \big[$. There exist $\rho_s \in ]0,\rho]$ and $C > 0$ which only depend on $s$ and such that, for $z \in \DD_+$, 
\begin{multline*} 
\big\| P(z) - P_0(z) \big\|_{\Lc(H_z^{s+1},H_z^{s-1})}\\
\leq C \left( \abs z^{2} \nr{G-\Id}_{\Sc^{-\bar \rho}} + \abs z^{2+\rho_s} \nr{aw}_{\Sc^{-1-\bar \rho}}  + \abs z^{2+\rho} \nr{w-1}_{\Sc^{-\bar \rho}} \right).
\end{multline*}
\end{proposition}

\begin{proof}
We apply Proposition \ref{prop:dec-sob} and Remark \ref{rem:dec-sob}. For $j,k \in \Ii 1 d$, we have by \eqref{hyp-swa} 
\[
\nr{\partial_j (G_{j,k}-\d_{j,k}) \partial_k}_{\Lc(H_z^{s+1},H_z^{s-1})} \leq \abs z^2 \nr{G_{j,k}-\d_{j,k}}_{\Lc(H_z^{s},H_z^{s})} \lesssim \abs z^2 \nr{G_{j,k}-\d_{j,k}}_{\Sc^{-\bar \rho}}.
\]
There exist $s_1,s_2 \in [s+1,s-1] \cap \big]-\frac d 2,\frac d 2\big[$, such that $s_1 = s_2 - \rho$. Then 
\[
\nr{z^2(w-1)}_{\Lc(H_z^{s+1},H_z^{s-1})} \leq \abs z^2 \nr{w-1}_{\Lc(H_z^{s_2},H_z^{s_1})} \lesssim \abs z^{2+\rho} \nr{w-1}_{\Sc^{-\bar \rho}}.
\]
Finally, we choose $\rho_s \in ]0,\rho]$ such that $[s-\rho_s,s + \rho_s] \subset \big]-\frac d 2,\frac d 2 \big[$. Then, we can find $s_1,s_2 \in [s-1,s+1] \cap \big]-\frac d 2,\frac d 2 \big[$, such that $s_2-s_1=1+\rho_s.$  By \eqref{hyp-swa}-\eqref{hyp-a}, we have
\begin{equation*}
\nr{zaw}_{\Lc(H_z^{s+1},H_z^{s-1})} \leq \abs z \nr{aw}_{\Lc(H_z^{s_2},H_z^{s_1})} \lesssim  \abs z^{2+\rho_s} \nr{aw}_{\Sc^{-1-\bar \rho}},
\end{equation*}
which gives the contribution of the damping term.
\end{proof}

In Proposition \ref{prop:reg-Ra} below, we will need smallness of $\abs z^{-2} (P(z)-P_0(z))$ in $\Lc(H_z^{s+1},H_z^{s-1})$ to get ellipitic regularity for $R(z)$. By Proposition \ref{prop:Pa-Po}, the contributions of $zaw$ and $z^2(w-1)$ are small for $z$ small. But this is not the case for the contribution of $(-\D_G + \D)$, unless $\nr{G-\Id}_{\Sc^{-\bar \rho}}$ is small. Since we have not assumed that this is the case, we will write the perturbation $G-\Id$ as a sum of a small perturbation and a compactly supported contribution, which will be handled differently. The following statement is Lemma 3.6 in \cite{Royer24} (notice that it is for this lemma that we cannot take $\bar \rho = \rho_0$).

\begin{lemma} \label{lem:chi-0}
Let $\gamma > 0$.  Then we can write $G = G_0 + G_\infty$ where $G_0 \in C_0^\infty$ and
$
\nr{G_\infty - \Id}_{\Sc^{-\bar \rho}} \leq \gamma.
$
\end{lemma}

Throughout the proofs of the two sections that follow, we will use commutators of the different operators involved with the operators of multiplication by the variables $x_j$ and the generator of dilations $\A$.\\

Let $T$ be a linear operator on $\Sc$. For $r \in ]0,1]$ and $j \in \Ii 1 d$, we set $\ad_{rx_j}(T) = T r x_j - r x_j T : \Sc \to \Sc$.  Then for $\mu = (\mu_1,\dots,\mu_d) \in \N^d$, we set
\[
\ad_{rx}^\mu = \ad_{r x_1}^{\mu_1} \circ \dots \circ \ad_{r x_d}^{\mu_d}.
\]
Notice that $\ad_{rx_j}$ and $\ad_{rx_k}$ commute for $j,k \in \Ii 1 d$.
For $z \in \DD_+$ and $j \in \Ii 1 d$ we set $\ad_{j,z} = \ad_{\abs z x_j}$. We also set $\ad_{0,z}(T) = \ad_{\A} (T) = T\A - \A T$. Finally, for $k \in \N$ and $J = (j_1,\dots,j_k) \in \Ii 0 d^k$, we set
\[
\ad_z^J(T) = \big(\ad_{j_1,z} \circ \dots \circ \ad_{j_k,z} \big) (T).
\]

Given $N \in \N$, we set $\Ic_N = \bigcup_{k=0}^N \Ii 0 d^k$. For $s_1,s_2 \in \R$, we say that $T$ belongs to $\Cc^N_z(H_z^{s_1},H_z^{s_2})$ if the operator $\ad_z^J(T)$ extends to a bounded operator from $H_z^{s_1}$ to $H_z^{s_2}$ for all $J \in \Ic_N$. In this case, we set
\[
\nr{T}_{\Cc^N_z(H_z^{s_1},H_z^{s_2})} = \sum_{J \in \Ic_N} \nr{\ad_z^J(T)}_{\Lc(H_z^{s_1},H_z^{s_2})}.
\]
We write $\nr{T}_{\Cc^N_z(H_z^{s})}$ for $\nr{T}_{\Cc^N_z(H_z^{s},H_z^{s})}$. Notice that for $s_1,s_2,s_3 \in \R$, $T_1 \in C^N_z(H_z^{s_1},H_z^{s_2})$ and $T_2 \in C^N_z(H_z^{s_2},H_z^{s_3})$, we have
\begin{equation} \label{eq:comm-T1T2}
\nr{T_2 T_1}_{\Cc^N_z(H_z^{s_1},H_z^{s_3})} \leq \nr{T_1}_{\Cc^N_z(H_z^{s_1},H_z^{s_2})}\nr{T_2}_{\Cc^N_z(H_z^{s_2},H_z^{s_3})}.
\end{equation}

We recall that the commutators of $A$ with derivatives and multiplication operators are given by
\begin{equation} \label{eq:comm-A-V-der}
[V,A] = i x \cdot \nabla V
\quad \text{and} \quad 
[\partial_j , \A] = -i  \partial_j.
\end{equation}
In particular, by induction on $k \in \N$,
\begin{equation} \label{eq:comm-Ak-xl}
\A^k x_j = x_j ( \A - i)^k.
\end{equation}

Now we can estimate the commutators with the main quantities involved in our study.

\begin{lemma} \label{lem:commutators}
Let $N \in \N$ and $s \in \R$. There exists $C > 0$ such that the following assertions hold for all $z \in \DD_+$.
\begin{enumerate}[\rm(i)]
\item If $s \in \big]-\frac d 2 , \frac d 2 \big[$, we have $\nr{G}_{\Cc^N_z(H_z^s)} \leq C$ and $\nr{w}_{\Cc^N_z(H_z^s)} \leq C$.
\item If $s \in \big]-\frac d 2 + \rho, \frac d 2 \big[$ then $\nr{G-\Id}_{\Cc^N_z(H_z^s,H_z^{s-\rho})} \leq C \abs z^{\rho}$ and  $\nr{w-1}_{\Cc^N_z(H_z^s,H_z^{s-\rho})} \leq C \abs z^{\rho}$.
\item Let $\rho' \in \{0,\rho\}$. If $s \in \big]-\frac d 2 + \rho'+1, \frac d 2 \big[$ then $\nr{a}_{\Cc^N_z(H_z^s,H_z^{s-\rho'-1})} \leq C \abs z^{1+\rho'}$.
\item For $j \in \Ii 1 d$ we have $\nr{\partial_j}_{\Cc^N_z(H_z^s,H_z^{s-1})} \leq  C \abs z$ and $\nr{\partial_j}_{\Cc^N_z(H_z^{s+1},H_z^{s})} \leq  C \abs z$.
\end{enumerate}
\end{lemma}

\begin{proof}
All these estimates are again given by Proposition \ref{prop:dec-sob} and Remark \ref{rem:dec-sob}. Notice that when $a = 0$, this is Lemma 3.7 in \cite{Royer24}, so we only prove the third statement. By \eqref{hyp-a} and \eqref{eq:comm-A-V-der}, we have
\begin{align*}
\nr{a}_{\Cc^N_z(H_z^s,H_z^{s-\rho'-1})}
 \lesssim \sum_{m=0}^N \nr{(x\cdot \nabla)^m a}_{\Lc(H_z^s,H_z^{s-1-\rho'})}
 \lesssim \abs z^{1+\rho'} \sum_{m=0}^N \nr{(x\cdot \nabla)^m (a)}_{\Sc^{-1-\bar \rho}}.  \qedhere
\end{align*}
\end{proof}

Now we can estimate the commutators with $P(z)$, $\theta_\s(z)$ (defined in \eqref{def:theta-sigma}) and $\gamma_j(z)$ (defined in \eqref{def:gamma01} and \eqref{def:gamma2}).

\begin{proposition} \label{prop:Pa-Po-comm}
Let $N \in \N$.
\begin{enumerate}[\rm(i)]
\item Let $s \in \big] - \frac d 2, \frac d 2 \big[$. There exists $C > 0$ such that, for $z \in \DD_+$, 
\[
\nr{P(z)}_{\Cc^N_z(H_z^{s+1},H_z^{s-1})} \leq C \abs z^{2}.
\]

\item Let $\sigma \in \{0,1,2\}$. We set $\s' = 0$ if $\s \in \{0,1\}$ and $\s' = 1$ if $\s = 2$. Let $s \in \big] - \frac d 2 -\s' + \s + \rho, \frac d 2+\s' \big[$. Then there exists $C > 0$ such that for $z \in \DD_+$ we have 
\[
\nr{\th_\sigma(z)}_{\Cc^N_z(H_z^{s},H_z^{s-\s-\rho})} \leq C \abs z^{\s + \rho}.
\]

\item Let $j \in \{0,1,2\}$ and $j' = \min(1,j) $. Then for $s \in \big] - \frac d 2 + j', \frac d 2 \big[$ there exists $C > 0$ such that for $z \in \DD_+$ we have
\begin{equation} \label{eq:estim-gamma-sigma}
\nr{\g_j(z)}_{\Cc^N_z(H_z^{s},H_z^{s-j'})} \leq C \abs z^{j}.
\end{equation}
\end{enumerate}
\end{proposition}

\begin{proof}
All these statements follow from Lemma \ref{lem:commutators} and \eqref{eq:comm-T1T2}. We detail for instance the estimate of $\th_2(z)$. In this case, we can take $s$ and $s-2-\rho$ in $\big]-\frac d 2 - 1,\frac d 2 +1\big[$. Notice that, compared to the proof of Proposition \ref{prop:Pa-Po}, here we can pay $2 + \rho$ derivatives and not only 2. Otherwise, the proof is similar. Let $s_1,s_2 \in [s-2-\rho,s] \cap \big]-\frac d 2,\frac d 2 \big[$ such that $s_2 - s_1 = 1+ \rho$. We have 
\begin{eqnarray*}
\lefteqn{\nr{\th_2(z)}_{\Cc^N_z(H_z^{s},H_z^{s-2-\rho})} = \nr{P(z)-P_0(z)}_{\Cc^N_z(H_z^{s},H_z^{s-2-\rho})}}\\
&& \lesssim \nr{\divg(G-\Id) \nabla }_{\Cc^N_z(H_z^{s},H_z^{s-2-\rho})} + \abs z \nr{aw}_{\Cc^N_z(H_z^{s},H_z^{s-2-\rho})} + \abs z^2 \nr{w-1}_{\Cc^N_z(H_z^{s},H_z^{s-2-\rho})} \\
&& \lesssim \nr{G-\Id}_{\Cc^N_z(H_z^{s-1},H_z^{s-1-\rho})} + \abs z \nr{aw}_{\Cc^N_z(H_z^{s_2},H_z^{s_1})} + \abs z^2 \nr{w-1}_{\Cc^N_z(H_z^{s-1},H_z^{s-1-\rho})} \\
&& \lesssim \abs{z}^{2+\rho} .
\end{eqnarray*}
The other estimates are similar. For the third statement, we observe that in $\g_2(z)$ there is a factor $\abs z$, which is why we only pay one derivative to get an estimate of order $\abs z^2$.
\end{proof}

\section{Elliptic regularity} \label{sec:elliptic-regularity}

In this section, we prove Propositions \ref{prop:regularity-Rz} and \ref{prop:regularity-Rz-A}. The parameter $\rho \in [0,\rho_0[$ is fixed by these statements, and we have also fixed $\bar \rho \in ]\rho,\rho_0[$. Proposition \ref{prop:regularity-Rz} is a consequence of Proposition \ref{prop:elliptic} below, and the proof of Proposition \ref{prop:regularity-Rz-A} is given at the end of the section. \\

Let $s \in \R$. For $r \in ]0,1]$ the resolvent $R_0(ir) = r^{-2} (D_r^2+1)\inv$ defines a bounded operator from $H_r^{s-1}$ to $H_r^{s+1}$ with norm $r^{-2}$. More generally, if we set
\[
\DI = \set{z \in \DD_+ \st  \mathsf{arg}(z) \in \left[\frac \pi 6 , \frac {5\pi}{6} \right]},
\]
then there exists $c_0 > 0$ such that for $s \in \R$ and $z \in \DI$, we have
\begin{equation} \label{eq:res-Po}
\nr{R_0(z)}_{\Lc(H_z^{s-1},H_z^{s+1})} \leq \frac {c_0} {\abs z^2}.
\end{equation}
Then, for $k \in \N^*$ and $s,s' \in \R$ such that $s' - s \leq k+2$ there exists $c > 0$ such that for $z \in \DI$, we have by \eqref{eq:der-Rz}-\eqref{eq:der-k-n} (for $R_0(z)$)
\begin{equation} \label{eq:res-Po-multiple}
\big\|R_0^{(k)}(z)\big\|_{\Lc(H_z^{s},H_z^{s'})} \leq  \frac{c}{\abs z^{k+2}}.
\end{equation}

Our purpose is to prove this kind of elliptic estimates for $R(z)$. By the usual elliptic regularity, this holds for any fixed $z \in \DD_+$ (see Proposition \ref{prop:Rz-base}). The difficulty is to get uniform estimates for $z$ close to 0.\\

We cannot extend \eqref{eq:res-Po} to $R(z)$ in full generality. We begin with the case $s = 0$.

\begin{proposition} \label{prop:Rari-H1}
There exists $c > 0$ such that for all $z \in \DI$, we have
\begin{equation*}
\nr{R(z)}_{\Lc(H_z\inv,H_z^1)} \leq \frac c {\abs z^2}.
\end{equation*}
More generally, for $N \in \N$ there exists $c_N > 0$ such that for $z \in \DI$, we have
\[
\nr{R(z)}_{\Cc^N_z(H_z\inv,H_z^1)} \leq \frac {c_N} {\abs z^2}.
\]
\end{proposition}

\begin{proof}
Let $z \in \DI$. The angle $\vartheta_z$ that appears in \eqref{eq:Lax-Milgram} now belongs to $\big[-\frac \pi 3,\frac \pi 3\big]$, so $\cos(\vartheta_z) \geq \frac 12$ for all $z \in \DI$. Then there exists $C>0$ such that for all $z \in \DI$ and $u \in H_z^1$, we have
\begin{align*}
\Re \innp{e^{-i\vartheta_z} P(z) u}{u}_{H\inv_z,H^1_z}
& = \cos(\vartheta_z)  \big(\innp{G \nabla u}{\nabla u}_{L^2} + \abs z^2 \innp{wu}{u}_{L^2} \big)+\abs z \innp{aw u}{u}_{L^2} \\
& \geq C \left( \nr{\nabla u}_{L^2}^2 + \abs z^2  \nr{u}_{L^2}^2    \right)
\\& \geq C \abs z^2 \nr{u}_{H^1_z}^2.
\end{align*}
The Lax-Milgram Theorem gives the first estimate.

Now let $J \in \Ic_N$. Then $\ad_z^J(R(z))$ can be written as a sum of terms of the form
\[
R(z) \ad_z^{J_1}(P(z)) R(z) \ad_z^{J_2}(P(z)) \dots R(z) \ad_z^{J_\nu}(P(z)) R(z),
\]
where $\nu \in \N$ and $J_1,\dots,J_\nu \in \Ic_N$. The second estimate then follows from Proposition \ref{prop:Pa-Po-comm}.
\end{proof}

We have a result similar to \eqref{eq:res-Po} if $G$ is a small perturbation of the flat metric (see the discussion after Proposition \ref{prop:Pa-Po}) and $s$ is not too large:

\begin{proposition} \label{prop:reg-Ra}
Let $s_1 \in \big]-\frac d 2,\frac d 2 +1 \big[$ and $s_2 \in \big]-\frac d 2-1,\frac d 2 \big[$ be such that $s_1 - s_2 \leq 2$. Let $N \in \N$. There exist $\gamma > 0$ and $c > 0$ such that if $\nr{G-\Id}_{\Sc^{-\bar \rho}} \leq \gamma$ then for $z \in \DI$, we have
\[
\nr{R(z)}_{\Cc_z^N(H_z^{s_2},H_z^{s_1})} \leq \frac c {\abs z^2}.
\]
\end{proposition}

\begin{proof}
Let $s_1' \in \big[s_1,\frac d 2 + 1 \big[$ and $s_2' \in \big]-\frac d 2 - 1 ,s_2 \big]$ be such that $s'_1 - s'_2 = 2$. It is enough to prove the estimate in $\Cc_z^N(H_z^{s_2'},H_z^{s_1'})$ instead of $\Cc_z^N(H_z^{s_2},H_z^{s_1})$.

Let $\rho_s > 0$ be given by Proposition \ref{prop:Pa-Po} applied with $s = \frac 12(s_1'+s_2')\in  \big]-\frac d 2,\frac d 2  \big[$. Then there exists $C > 0$ such that for $z \in \DI$, we have
\begin{align*}
     \nr{P(z)-P_0(z)}_{\Lc(H_z^{s'_1},H_z^{s'_2})} & =\nr{P(z)-P_0(z)}_{\Lc(H_z^{s+1},H_z^{s-1})}
    \\&\leq  C \abs z^2 \big(\gamma + \abs z^{\rho_s} \nr{aw}_{\Sc^{-\tilde \rho -1}} +\abs z^{\rho} \nr{w-1}_{\Sc^{-\tilde \rho}} \big).
\end{align*}
Thus, if $\gamma$ and $\abs z$ are small enough, we have
\begin{align*}
    \nr{P(z)-P_0(z)}_{\Lc(H_z^{s'_1},H_z^{s'_2})} \leq \frac {\abs z^2}{2c_0},
\end{align*}
where $c_0>0$ is given by \eqref{eq:res-Po}. This proves that 
\[
P(z) = P_0(z) \big( \Id + R_0(z) (P(z) - P_0(z)) \big) \quad  \in \Lc(H_z^{s'_1},H_z^{s'_2})
\]
is boundedly invertible for $z$ small enough, and $R(z)$ can be seen as an operator in $\Lc(H_z^{s'_2},H_z^{s'_1})$ with
\[
\nr{R(z)}_{\Lc(H_z^{s'_2},H_z^{s'_1})} \leq \frac {2 c_0}{\abs z^2}.
\]
For $z \in \DI$ outside a neighborhood of 0, the same estimate holds (with a possibly different constant) by elliptic regularity. Thus, the estimate is proved for $N = 0$. We estimate the commutators as in Proposition \ref{prop:Rari-H1}, and the proof is complete.
\end{proof}

For $k \in \N$ and $\jj \in \{0,1,2\}^k$, we set
\[
m'(k,\jj) =  2(k+1)- \sum_{\ell = 1}^k \min(1,j_\ell).
\]
In particular,
\[
m'(k,\jj) \geq m(k,\jj).
\]
The motivation for this notation is that $m'(k,\jj)$ measures the gain of regularity for the operator $\Rc_{k,\jj}(z)$ when $z \in \DI$.

\begin{proposition} \label{prop:reg-T}
Let $k \in \N$ and $ \jj=(j_1,\dots,j_k) \in \{0,1,2\}^k$. Let $s_1 \in \big]-\frac d 2,\frac d 2 +1 \big[$ and $s_2 \in \big]-\frac d 2-1,\frac d 2 \big[$ be such that
\[
s_1 - s_2 \leq m'(k,\jj).
\]
Let $N \in \N$. Let $\gamma > 0$ be given by Proposition \ref{prop:reg-Ra}. There exists $C > 0$ such that if $\nr{G-\Id}_{\Sc^{-\bar \rho}} \leq \gamma$ then for $z \in \DD_+$ and $z' \in \DI$ with $\abs z = \abs {z'}$,
\[
\nr{\Rc_{k,\jj}(z',z)}_{\Cc_z^N(H_z^{s_2},H_z^{s_1})} \leq C \abs z^{-m(k,\jj)}.
\]
\end{proposition}

\begin{proof}
Let $\mathfrak s \in ]0,\frac d 2[$ such that $\mathfrak s \geq s_1 - 1$. For $\ell \in \{1,\dots,k\}$ we set $j'_\ell = \min(1,j_\ell)$. We set $\s_k = \min(\mathfrak s, s_2 + 2)$. By Proposition \ref{prop:reg-Ra}, we have
\[
\nr{R(z')}_{\Cc_z^N(H_z^{s_2},H_z^{\s_k})} \lesssim \abs z^{-2}.
\]
We define $\s_{k-1},\s_{k-2},\dots,\s_2,\s_1$ inductively by
\[
\s_{\nu} = \min(\mathfrak s,\s_{\nu+1} - j'_{\nu+1} + 2), \quad \nu = k-1,\dots,1.
\]
In particular,
\[
\frac d 2 > \s_1 \geq \dots \geq \s_{k} > -\frac d 2 + 1.
\]
By Propositions \ref{prop:Pa-Po-comm} and \ref{prop:reg-Ra}, we have for all $\nu \in \Ii 2 {k-1}$
\[
\nr{R(z')\g_{j_\nu}(z)}_{\Cc_z^N(H_z^{\s_{\nu+1}},H_z^{\s_\nu})} \lesssim \abs z^{j_\nu-2}.
\]
Finally, if $\s_1 \neq \mathfrak s$ then $\s_\ell\neq \mathfrak s$ for all $\ell \in \{1,\dots,k\}$ and
\begin{align*}
\s_1 = s_2 + 2  + \sum_{\nu = 1}^{k-1} (2-j'_{\nu+1}) = s_2 + m'(k,\jj) -(2-j'_1).
\end{align*}
In any case we have $s_1-\s_1 \leq 2 - j'_1$, so
\[
\nr{R(z')\g_{j_1}(z)}_{\Cc_z^N(H_z^{\s_1},H_z^{s_1})} \lesssim \abs z^{j_1-2}.
\]
All these estimates together give
\[
\nr{\Rc_{k,\jj}(z',z)}_{\Cc_z^N(H_z^{s_2},H_z^{s_1})} \lesssim \abs z^{(j_1+\dots+j_k) - 2(k+1)} = \abs z^{-m(k,\jj)},
\]
and the proof is complete.
\end{proof}

Notice in this proof that if $s_1 \geq \frac d 2 -  j_1'$ (or $s_2 \leq -\frac d 2 + j_k'$) it is important that the first (last) factor in $\Rc_{k,\jj}(z',z)$ is a resolvent $R(z')$.

The results of Proposition \ref{prop:reg-T} also hold in the particular case $G=\Id$, $a=0$ and $w=1$, that is, for $\Rc^0_{k,\jj}(z',z)$ instead of $\Rc_{k,\jj}(z',z)$.

Now we consider the case where we have resolvents $R(z')$ and $R_0(z')$ with an inserted factor $\th_\sigma(z)$. For further use, we also consider the possibility of having a derivative next to $\Rc_{k,\jj}(z',z)$.

\begin{proposition} \label{prop:reg-T-T0}
Let $k,k_0 \in \N$, $ \jj \in \{0,1,2\}^k$ and $ \jj_0 \in \{0,1,2\}^{k_0}$.
Let  $\a \in \N^d$ with $\abs{\a} \leq 1$. Let $s_1,s_2 \in \big]-\frac d 2,\frac d 2 \big[$
 such that
\begin{equation} \label{eq:hyp-s1s2}
s_1 - s_2 \leq -\abs{\a}+ m'(k,\jj) - (\s+\rho) + m'(k_0,\jj_0).
\end{equation}
Let $N \in \N$. Let $\gamma > 0$ be given by Proposition \ref{prop:reg-Ra}. Assume that $\a = 0$ or $s_1 < \frac d 2 - \rho$. There exists $C > 0$ such that if $\nr{G-\Id}_{\Sc^{-\bar \rho}} \leq \gamma$ then for $z \in \DD_+$ and $z' \in \DI$ with $\abs z = \abs {z'}$, we have
\[
\nr{D^{\a} \Rc_{k,\jj}(z',z) \th_\sigma(z) \Rc^0_{k_0,\jj_0}(z',z)}_{\Cc_z^N(H_z^{s_2},H_z^{s_1})} \leq C \abs z^{-m(k,\jj) + (\s+\rho) - m(k_0,\jj_0)+\abs{\a}}.
\]
\end{proposition}

\begin{proof}
As in Proposition \ref{prop:Pa-Po-comm}, we set $\s' = 1$ if $\s = 2$ and $\s'=0$ if $\s \in \{0,1\}$. In particular, $\s \leq \s'+1$. We recall that $m'(k,\jj), m'(k_0,\jj_0) \geq 2$. Since $\a = 0$ or $s_1 < \frac d 2 - \rho$, we have
\[
s_1 + \abs \a -m'(k,\jj) + \s + \rho < \frac d 2 + \s'.
\]
On the other hand, we also have
\[
s_2+ m'(k_0,\jj_0) > -\frac d 2 - \s' + \s + \rho.
\]
Then we can consider
\[
\mathfrak s_+ \in \left]-\frac d 2 - \s' + \s + \rho , \frac d 2 + \s'\right[ \cap \big[ s_1 + \abs \a -m'(k,\jj) + \s + \rho, s_2 + m'(k_0,\jj_0) \big]
\]
(notice that these two intervals are not empty since $d + 2\s'-\s-\rho > 0$ and by \eqref{eq:hyp-s1s2}), and we set
\[
\mathfrak s_- = \mathfrak s_+ - (\s + \rho).
\]
Since $\mathfrak s_+ - s_2 \leq m'(k_0,\jj_0)$ we have, by Proposition \ref{prop:reg-T} applied to $\Rc^0_{k_0,\jj_0}(z',z)$,
\[
\nr{ \Rc^0_{k_0,\jj_0}(z',z)}_{\Cc_z^N(H_z^{s_2},H_z^{\mathfrak s_+})} \lesssim \abs z^{-m(k_0,\jj_0)}.
\]
By Proposition \ref{prop:Pa-Po-comm}, we have
\[
\nr{\th_\sigma(z)}_{\Cc_z^N(H_z^{\mathfrak s_+},H_z^{\mathfrak s_-})} \lesssim \abs z^{\sigma + \rho}.
\]
Since $s_1 + \abs \a - \mathfrak s_- \leq m'(k,\jj)$ we have, by Proposition \ref{prop:reg-T} again,
\[
\nr{\Rc_{k,\jj}(z',z)}_{\Cc_z^N(H_z^{\mathfrak s_-},H_z^{s_1+\abs \a})} \lesssim \abs z^{-m(k,\jj)}.
\]
And finally, $\nr{D^\a}_{\Cc_z^N(H_z^{s_1+ \abs \a},H_z^{s_1})} \lesssim \abs z^{\abs \a}$. These four estimates together give the result.
\end{proof}

The previous results hold only if $G$ is close to $\Id$, in the sense that $\nr{G-\Id}_{\Sc^{-\bar \rho}} \leq \g$, which is not the case in general. The idea will be to prove first all the estimates for a small perturbation of the flat metric and then to add the contribution of a compactly supported perturbation with one of the following two decompositions for $R_{k,\jj}(z',z)$.

\begin{lemma} \label{lem:R-Rinf}
Let $\gamma > 0$ be given by Proposition \ref{prop:reg-Ra}. Let $G_\infty$ and $G_0$ be given by Lemma \ref{lem:chi-0}. Let $z \in \DD_+$ and $z' \in \DI$. Let $k \in \N$ and $ \jj = (j_1,\dots,j_k) \in \{0,1,2\}^k$. We define $R_\infty(z)$ and $\Rc_{k,\jj}^\infty(z',z)$ as $R(z)$ and $\Rc_{k,\jj}(z',z)$ with $G$ replaced by $G_\infty$.
\begin{enumerate} [\rm(i)]
\item $\Rc_{k,\jj}(z',z)$ can be written as a sum of terms of the form
\begin{equation} \label{eq:terme-T}
\Rc_{\k_0,\ll_0}^\infty(z',z) B_1(z') \Rc_{\k_1,\ll_1}^\infty(z',z) B_2(z') \dots B_{\nu}(z')\Rc_{\k_\nu,\ll_\nu}^\infty(z',z),
\end{equation}
where $\nu \in \N$, $B_1(z'),\dots,B_\nu(z')$ belong to $\{\D_{G_0},\D_{G_0} R(z')\D_{G_0}\}$, and $\k_0,\dots,\k_\nu \in \N$ and $\ll_i \in \{0,1,2\}^{\k_i}$, $i \in \Ii 0 \nu$, are such that
\begin{equation} \label{eq:mu-muj}
 \sum_{i = 0}^{\nu}  m(\k_i,\ll_i)- 2\nu  = m(k,\jj).
\end{equation}
Moreover, an operator of the form $D_{\ell_1} R(z') D_{\ell_2}$ for $1 \leq \ell_1,\ell_2 \leq d$ extends to a bounded operator on $L^2$ uniformly in $z' \in \DI$.
\item We have
\begin{equation} \label{eq:R-Rinf-2}
\Rc_{k,\jj}(z',z)
= \Rc_{k,\jj}^\infty(z',z) + \sum_{\nu=0}^k \Rc_{\nu,\ll_\nu}(z',z) \D_{G_0} \Rc_{k-\nu,\ll_\nu^*}^\infty(z',z),
\end{equation}
where for $\nu \in \Ii 0 k$ we have set $\ll_\nu=(j_1,\dots,j_\nu)$ and $\ll_\nu^* = (j_{\nu+1},\dots,j_k)$, and we have
\begin{equation} \label{eq:m-nu}
m(\nu,\ll_\nu) + m(k-\nu,\ll_\nu^*) = m(k,\jj) + 2.
\end{equation}
\end{enumerate}
\end{lemma}

\begin{proof}
We obtain the first statement by iterating the resolvent identity
\begin{equation*}
R(z') = R_\infty(z') + R_\infty(z') \D_{G_0} R_\infty(z') + R_\infty(z') \D_{G_0} R(z') \D_{G_0} R_\infty(z'),
\end{equation*}
so that all the remaining factors $R(z')$ are hidden in the factors $B_\ell(z')$. We have \eqref{eq:mu-muj} since each time we get a factor $B_\ell(z')$ we increase by one the number of resolvents (roughly, in the corresponding terms we replace one resolvent $R(z')$ by two resolvents $R_\infty(z')$). The boundedness of $D_{\ell_1} R(z') D_{\ell_2}$ comes from Proposition \ref{prop:Rari-H1} and \eqref{eq:derivee-alpha}.

The equality \eqref{eq:R-Rinf-2} simply follows from the resolvent identity $R(z') = R_\infty(z') + R(z') \D_{G_0} R_\infty(z')$.
\end{proof}

In the following proofs, we will also use the following simple remark. For $p \in \N^*$ and $\eta_0,\eta_1,\dots,\eta_p \geq 0$ we have
\begin{equation} \label{eq:min}
\min(\eta_0 - \eta_1,0) + \dots + \min(\eta_0 - \eta_p,0) \geq \min\big(\eta_0-(\eta_1+\dots + \eta_p),0\big).
\end{equation}

Now we can prove the estimates for $\Rc_{k,\jj}(z',z)$ and $\Rc_{k,\jj}(z',z) \th_\sigma(z)  \Rc^0_{k_0,\jj_0}(z',z)$ without the smallness assumption for $G-\Id$, which gives Proposition \ref{prop:regularity-Rz}. The derivatives that we add in \eqref{eq:R-derivees} will be useful for the proof of \eqref{eq:R-th-R0}.

\begin{proposition} \label{prop:elliptic}
Let $s_1,s_2 \in \big[0,\frac d 2\big[$, $\d_1 > s_1$ and $\d_2 > s_2$. Let $k,k_0 \in \N$, $\jj \in \{0,1,2\}^k$ and $\jj_0 \in \{0,1,2\}^{k_0}$. Let $\sigma \in \{0,1,2\}$. Let $\a_1,\a_2 \in \N^d$ with $\abs{\a_1},\abs{\a_2} \leq 1$.
There exists $C > 0$ such that for $z \in \DD_+$ and $z' \in \DI$ with $\abs z = \abs {z'}$, we have
\begin{equation} \label{eq:R-derivees}
\nr{\pppg x^{-\d_1} D^{\a_1} \Rc_{k,\jj}(z',z) D^{\a_2} \pppg x^{-\d_2}}_{\Lc(L^2)} \leq C \abs {z}^{\min(s_1+s_2-m(k,\jj)  + \abs {\a_1} + \abs {\a_2},0)}
\end{equation}
and
\begin{equation} \label{eq:R-th-R0}
 \nr{\pppg x^{-\d_1}   \Rc_{k,\jj}(z',z) \th_\sigma(z)  \Rc^0_{k,\jj_0}(z',z) \pppg x^{-\d_2}}_{\Lc(L^2)}
\leq C \abs z^{\min(s_1 + s_2  -m(k,\jj) + (\s + \rho) -m(k_0,\jj_0),0)}.
\end{equation}
\end{proposition}

\begin{proof}
Let $\gamma > 0$ be given by Proposition \ref{prop:reg-Ra}. We use the notation of Lemma \ref{lem:R-Rinf}. Then Propositions \ref{prop:reg-T} and \ref{prop:reg-T-T0} apply to $R_\infty(z')$. \\
\stepp We assume that $s_1 + s_2 \leq m(k,\jj) - \abs {\a_1} - \abs {\a_2}$. Otherwise we can always replace $s_1$ and $s_2$ by $\tilde s_1 \in [0,s_1]$ and $\tilde s_2 \in [0,s_2]$ such that $\tilde s_1 + \tilde s_2 = m(k,\jj)  - \abs {\a_1} - \abs {\a_2}$. \\
By Proposition \ref{prop:reg-T}, \eqref{eq:derivee-alpha} and Lemma \ref{lem:weight}, the estimate \eqref{eq:R-derivees} holds with $\Rc_{k,\jj}(z',z)$ replaced by $ \Rc_{k,\jj}^\infty(z',z)$.

To prove \eqref{eq:R-derivees} in general, we estimate a term $\Rc(z',z)$ of the form \eqref{eq:terme-T}. Using \eqref{eq:R-derivees} proved for $R_\infty(z')$, the compactness of the support of $G_0$, the derivatives given by the operator $\D_{G_0}$ and the boundedness of an operator of the form $\partial_{x_{j_1}} R(z') \partial_{x_{j_2}}$, we obtain
\[
\nr{\pppg x^{-\d_1} D^{\a_1} \mathsf \Rc(z',z) D^{\a_2} \pppg x^{-\d_2}}_{\Lc(L^2)} \lesssim \sum_{\ell_1,\dots,\ell_{2\nu} = 1}^d \Nc_{\ell_1,\dots,\ell_{2\nu}},
\]
where
\begin{eqnarray*}
\lefteqn{\Nc_{\ell_1,\dots,\ell_{2\nu}}}\\
&& \lesssim   \nr{\pppg x^{-\d_1} D^{\a_1} \Rc_{\k_0,\ll_0}^\infty(z',z) D_{\ell_1} \pppg x^{-\d_2}}  \times \prod_{i=1}^{\nu-1} \nr{\pppg x^{-\d_1} D_{\ell_{2i}} \Rc_{\k_i,\ll_i}^\infty(z',z)  D_{\ell_{2i+1}} \pppg x^{-\d_2}} \\
&& \quad \times
\nr{\pppg x^{- \d_1} D_{\ell_{2\nu}} \Rc_{\k_\nu,\ll_\nu}^\infty(z',z) D^{\a_2} \pppg x^{-\d_2}} \\
&& \lesssim  \abs z^{\min(s_1+s_2 + \abs{\a_1}+1-m(\k_0,\ll_0), 0)} \times \prod_{i=1}^{\nu-1}  \abs z^{\min(s_1+s_2 + 2 - m(\k_i,\ll_i), 0)} \\
&& \quad \times  \abs z^{\min(s_1+s_2 + \abs{\a_2}+1-m(\k_\nu,\ll_\nu), 0)}\\
&& \lesssim \abs z^{\min(s_1+s_2 + \abs{\a_1} + \abs{\a_2} + 2 \nu - \sum_{i=0}^\nu  m(\k_i,\ll_i), 0)}.
\end{eqnarray*}
For the last inequality, we have used \eqref{eq:min} with $\eta_0 = s_1 + s_2$, $\eta_1 = m(\k_0,\ll_0) -\abs{\a_1} - 1$, $\eta_j = m(\k_{j-1},\ll_{j-1}) - 2$ for $j =2,\dots,\nu$, and finally $\eta_{\nu+1} = m(\k_\nu,\ll_\nu) - \abs {\a_2}-1$. With \eqref{eq:mu-muj}, this proves \eqref{eq:R-derivees}.

\stepp
We turn to the proof of \eqref{eq:R-th-R0}. Let $\a \in \N^d$ with $\abs \a \leq 1$. If $\a = 0$ or $s_1 < \frac d 2 - \rho$, then by Proposition \ref{prop:reg-T-T0} and Lemma \ref{lem:weight}, we have
\begin{multline} \label{eq:R-th-R0-derivee}
\nr{\pppg x^{-\d_1}  D^{\a} \Rc_{k,\jj}^\infty(z',z) \th_\sigma(z)  \Rc^0_{k_0,\jj_0}(z',z) \pppg x^{-\d_2}}_{\Lc(L^2)}\\
\lesssim \abs z^{\min(s_1 + s_2 +\abs{\a} -m(k,\jj) + (\s+\rho)-m(k_0,\jj_0),0)}.
\end{multline}
Now that \eqref{eq:R-derivees} is proved, we no longer need \eqref{eq:terme-T}, and we use the simpler decomposition \eqref{eq:R-Rinf-2} instead. We have just proved \eqref{eq:R-th-R0} with $\Rc_{k,\jj}(z',z)$ replaced by $\Rc_{k,\jj}^\infty(z',z)$. We consider \eqref{eq:R-th-R0} with $\Rc_{k,\jj}(z',z)$ replaced by a term of the sum in \eqref{eq:R-Rinf-2}.
We set $s_1' = {\max(s_1 - \rho,0)}< \frac d 2 - \rho$ and consider $\d_1' > s_1'$. For $\nu \in \{0,\dots,k\}$ we have, by \eqref{eq:R-derivees} and \eqref{eq:R-th-R0-derivee},
\begin{eqnarray*}
\lefteqn{ \nr{\pppg x^{-\d_1}   \Rc_{\nu,\ll_\nu}(z',z) \D_{G_0} \Rc_{k-\nu,\ll_\nu^*}^\infty(z',z) \th_\sigma(z)  \Rc^0_{k_0,\jj_0}(z',z) \pppg x^{-\d_2}}}\\
&& \lesssim \sum_{\abs {\a_1},\abs {\a_2} = 1} \nr{\pppg x^{-\d_1}   \Rc_{\nu,\ll_\nu}(z',z) D^{\a_1} \pppg {x}^{-\d_2}} \\
&& \qquad \qquad \times \nr{\pppg x^{-\d_1'} D^{\a_2} \Rc_{k-\nu,\ll_\nu^*}^\infty(z',z) \th_\sigma(z)  \Rc^0_{k_0,\jj_0}(z',z) \pppg x^{-\d_2}}\\
&& \lesssim \abs z^{\min(s_1+s_2 -m(\nu,\ll_\nu) + 1,0)} \abs z^{\min(s_1'+s_2 +1 -m(k-\nu,\ll_\nu^*) + (\s+\rho) - m(k_0,\jj_0),0)}\\
&& \lesssim  \abs z^{\min(s_1 + s_2 -m(k,\jj) + (\s+\rho)-m(k_0,\jj_0),0)}.
\end{eqnarray*}
Here, we have applied \eqref{eq:min} with $\eta_0 = s_1' + s_2 \geq 0$, $\eta_1 = m(\nu,\ell_\nu) - 1 + s_1'-s_1 \geq 0$ and $\eta_2 = m(k-\nu,\ll_\nu^*) + m(k_0,\jj_0) - 1 - (\s+\rho) \geq 0$. Then we have used \eqref{eq:m-nu}. This concludes the proof.
\end{proof}

Now we turn to the proof of Proposition \ref{prop:regularity-Rz-A}. We first prove a statement with $\A$ replaced by $\abs z x$ and without the factors $\g_{\s'}(z)$. This will be combined with Proposition \ref{prop:xTA} below.

\begin{proposition} \label{prop:elliptic-2}
Let $s \in \big[0,\frac d 2\big[$ and $\d > s$. Let $k,k_0 \in \N$, $\jj \in \{0,1,2\}^k$ and $\jj_0 \in \{0,1,2\}^{k_0}$.
There exists $C > 0$ such that, for $z \in \DD_+$ and $r=\abs z$,
\begin{equation}
\label{eq:x-R-rx}
\nr{\pppg x^{-\d}  \Rc_{k,\jj}(ir,z) \pppg{rx}^\d}_{\Lc(L^2)} \leq C r^{\min(s-m(k,\jj),0)},
\end{equation}
\begin{equation}
\label{eq:x-R-th-R0-rx}
\nr{\pppg x^{-\d} \Rc_{k,\jj}(ir,z) \th_\sigma(z)
\Rc^0_{k_0,\jj_0}(ir,z) \pppg{rx}^\d}_{\Lc(L^2)} \leq C
r^{\min(s  -m(k,\jj)+ (\s +\rho) - m(k_0,\jj_0),0)},
\end{equation}
\begin{equation}
\label{eq:rx-R0-x}
\nr{\pppg{rx}^\d \Rc^0_{k_0,\jj_0}(ir,z)\pppg x^{-\d}}_{\Lc(L^2)} \leq  C r^{\min(s-m(k_0,\jj_0),0)},
\end{equation}
\begin{equation}
\label{eq:rx-R-th-R0-x}
\nr{\pppg{rx}^\d     \Rc_{k,\jj}(ir,z) \th_\sigma(z) \Rc^0_{k_0,\jj_0}(ir,z) \pppg x^{-\d}}_{\Lc(L^2)}  \leq   C r^{\min(s-m(k,\jj)+(\s+ \rho)-m(k_0,\jj_0)+(\s+ \rho),0)}.
\end{equation}
\end{proposition}

\dimm{On ne met pas de $\g$ sur les côtés, car ils seront à la proposition \ref{prop:xTA}. Si on veut considérer la dimension 2, on peut avoir à les mettre en partie de chaque côté}

\begin{proof}
As for the proof of Proposition \ref{prop:elliptic}, we consider $\gamma > 0$ given by Proposition \ref{prop:reg-Ra} and we use the notation of Lemma \ref{lem:R-Rinf}. We first prove the estimates with $\Rc_{k,\jj}(ir,z)$ replaced by $\Rc_{k,\jj}^{\infty}(ir,z)$ and with an additional derivative. Then we will deduce the general case with Lemma \ref{lem:R-Rinf}.

\stepp We begin with \eqref{eq:x-R-rx}. Let $\a \in \N^d$ with $\abs \a \leq 1$. Let $\mu \in \N$  and $\b \in \N^d$ with $\abs \b \leq 2 \mu$. We can write $\pppg {rx}^{-2\mu} D^{\a}\Rc_{k,\jj}^{\infty}(ir,z) (rx)^\b$ as a sum of terms of the form
\[
\pppg {rx}^{-2\mu} (rx)^{\b_1} \ad_{rx}^{\b_2} \big(D^{\a} \Rc_{k,\jj}^{\infty}(ir,z)\big),
\]
 where $\b_1 + \b_2 = \b$. We estimate such a term. Let $\mathfrak s = \min(s, m'(k,\jj)- \abs {\a})$. By Lemma \ref{lem:commutators}, Proposition \ref{prop:reg-T} and \eqref{eq:comm-T1T2}, we have
 \[
 \nr{\ad_{rx}^{\b_2} \big(D^{\a} \Rc_{k,\jj}^{\infty}(ir,z)\big)}_{\Lc(L^2,H_z^{\mathfrak s})} \lesssim r^{\abs \a-m(k,\jj)}.
 \]
Since $\pppg {rx}^{-2\mu} (rx)^{\b_1}$ is uniformly bounded in $\Lc(H_r^{\mathfrak s})$, we get
\[
\nr{\pppg {rx}^{-2\mu} (rx)^{\b_1} \ad_{rx}^{\b_2} \big(D^{\a} \Rc_{k,\jj}^{\infty}(ir,z)\big)}_{\Lc(L^2,H_r^{\mathfrak s})} \lesssim r^{\abs \a-m(k,\jj)}.
\]
This proves that for any $\mu \in \N$, we have
\[
\nr{\pppg {rx}^{-2\mu} D^{\a}\Rc_{k,\jj}^{\infty}(ir,z)  \pppg {rx}^{2\mu}}_{\Lc(L^2,H_r^\mathfrak{s})} \lesssim r^{\abs \a-m(k,\jj)}.
\]
By interpolation,
\begin{equation} \label{eq:rxd-rxd}
\nr{\pppg {rx}^{-\d} D^{\a} \Rc_{k,\jj}^{\infty}(ir,z) \pppg{rx}^\d}_{\Lc(L^2,H_r^{\mathfrak s})} \lesssim r^{\abs \a-m(k,\jj)}.
\end{equation}
On the other hand, by Lemma \ref{lem:weight},
\begin{equation} \label{eq:xd-rxd}
\begin{aligned}
\nr{\pppg x^{-\d} \pppg {rx}^\d}_{\Lc(H^\mathfrak{s}_r,L^2)}
& \lesssim
 \big\| (1+\abs{rx}^\d) \pppg{x}^{-\d} \big\|_{\Lc(H_r^\mathfrak{s},L^2)} \\
& \lesssim  \|\pppg {x}^{-\d}\|_{\Lc(H^\mathfrak{s}_r,L^2)} + r^\d \big\|\abs x^\d \pppg{x}^{-\d} \big\|_{\Lc(L^2)}\\
& \lesssim r^\mathfrak{s}.
\end{aligned}
\end{equation}
Finally,
\begin{eqnarray} \label{eq:x-R-rx-der}
\lefteqn{\nr{\pppg x^{-\d}  D^{\a} \Rc_{k,\jj}^{\infty}(ir,z)  \pppg{rx}^\d}_{\Lc(L^2)} }
\\
\nonumber &&\leq \nr{\pppg x^{-\d} \pppg {rx}^\d}_{\Lc(H^\mathfrak{s}_r,L^2)}\nr{\pppg {rx}^{-\d} D^{\a} \Rc_{k,\jj}^{\infty}(ir,z) \pppg{rx}^\d}_{\Lc(L^2,H_r^\mathfrak{s})} ,    \\
\nonumber && \lesssim r^{ \min(s+\abs \a-m(k,\jj),0)}.
\end{eqnarray}

Now we use \eqref{eq:R-Rinf-2} to prove \eqref{eq:x-R-rx}. The contribution of $\Rc_{k,\jj}^\infty(ir,z)$ is given by \eqref{eq:x-R-rx-der} applied with $\a = 0$. Then, for $\nu \in \{0,\dots,k\}$, we have by \eqref{eq:R-derivees} and \eqref{eq:x-R-rx-der}
\begin{eqnarray*}
\lefteqn{\nr{\pppg x^{-\d}  \Rc_{\nu,\ll_\nu}(z',z) \D_{G_0} \Rc_{k-\nu,\ll_\nu^*}^\infty(z',z) \pppg{rx}^\d}}\\
&& \lesssim \sum_{\abs{\a_1},\abs{\a_2} = 1} \nr{\pppg x^{-\d} \Rc_{\nu,\ll_\nu}(z',z) D^{\a_1}} \nr{\pppg x^{-\d} D^{\a_2} \Rc_{k-\nu,\ll_\nu^*}^\infty(z',z) \pppg{rx}^\d}\\
&& \lesssim r^{\min(s + 1 -m(\nu,\ll_\nu),0)}  r^{\min(s+1 - m(k-\nu,\ll_\nu^*),0)}\\
&& \lesssim r^{\min(s-m(k,\jj),0)}.
\end{eqnarray*}
Here we have applied \eqref{eq:min} with $\eta_0 = s$. This gives \eqref{eq:x-R-rx}. Similarly, we can prove
\begin{equation} \label{eq:rxd-R-xd}
\nr{ \pppg{rx}^\d \Rc_{k,\jj}(ir,z)  \pppg x^{-\d}}_{\Lc(L^2)} \leq C r^{\min(s-m(k,\jj),0)},
\end{equation}
which gives \eqref{eq:rx-R0-x} as a particular case.

\stepp We proceed similarly for \eqref{eq:x-R-th-R0-rx}. If $\a = 0$ or $s < \frac d 2 - \rho$, we get with Proposition \ref{prop:reg-T-T0}
\begin{equation} \label{eq:x-R-th-R0-rx-derivee}
\nr{\pppg x^{-\d}D^{\a} \Rc_{k,\jj}^{\infty}(ir,z) \th_\sigma(z) \Rc^0_{k_0,\jj_0}(ir,z)\pppg{rx}^\d}_{\Lc(L^2)} \lesssim
r^{\min(s + \abs \a  -m(k,\jj) + (\s + \rho) -m(k_0,\jj_0),0)}.
\end{equation}
We set $s' = \max(s-\rho,0)< \frac d 2 - \rho$ and consider $\d' > s'$. Then we get \eqref{eq:x-R-th-R0-rx} with Lemma \ref{lem:R-Rinf} since, for $\nu \in \{0,\dots,k\}$, we have by \eqref{eq:R-derivees} and \eqref{eq:x-R-th-R0-rx-derivee}
\begin{eqnarray*}
\lefteqn{\nr{\pppg x^{-\d}  \Rc_{\nu,\ll_\nu}(z',z) \D_{G_0} \Rc_{k-\nu,\ll_\nu^*}^\infty(z',z)\th_\sigma(z) \Rc^0_{k_0,\jj_0}(ir,z)  \pppg{rx}^\d}}\\
&& \lesssim \sum_{\abs{\a_1},\abs{\a_2} = 1} \nr{\pppg x^{-\d} \Rc_{\nu,\ll_\nu}(z',z) D^{\a_1}} \nr{\pppg x^{-\d'} D^{\a_2} \Rc_{k-\nu,\ll_\nu^*}^\infty(z',z) \th_\sigma(z) \Rc^0_{k_0,\jj_0}(ir,z)  \pppg{rx}^\d}\\
&& \lesssim r^{\min(s + 1 -m(\nu,\ll_\nu),0)}  r^{\min(s'+1 - m(k-\nu,\ll_\nu^*)+(\s+\rho) - m(k_0,\jj_0),0)}\\
&& \lesssim r^{\min(s-m(k,\jj) + (\s+\rho) - m(k_0,\jj_0),0)}.
\end{eqnarray*}
We have applied \eqref{eq:min} with $\eta_0 = s'$.

\stepp We finish with \eqref{eq:rx-R-th-R0-x}. We proceed again similarly. Let
\[
\mathfrak s = \min(s, m(k,\jj) - (\s+\rho) + m(k_0,\jj_0) - \abs \a).
\]
With Proposition \ref{prop:reg-T-T0} and the same strategy as for \eqref{eq:rxd-rxd} we get
\begin{equation*}
\nr{\pppg{rx}^\d  D^{\a} \Rc_{k,\jj}^{\infty}(ir,z) \th_{\sigma}(z) \Rc^0_{k_0,\jj_0}(ir,z) \pppg {rx}^{-\d}}_{\Lc(H^{-\mathfrak s}_r, L^2)}   \lesssim r^{\abs{\a}-m(k,\jj)  + (\s+\rho) -m(k_0,\jj_0)}.
\end{equation*}
As in \eqref{eq:xd-rxd}, or by duality, we obtain  $\nr{\pppg{rx}^\d \pppg x^{-\d}}_{\Lc(L^2,H_r^{-\mathfrak s})} \lesssim r^{\mathfrak s}$, and hence
\begin{equation*}
\nr{\pppg{rx}^\d   D^{\a}\Rc_{k,\jj}^{\infty}(ir,z) \th_{\sigma}(z) \Rc^0_{k_0,\jj_0}(ir,z) \pppg x^{-\d}}_{\Lc(L^2)}  \lesssim r^{\min(s+\abs{\a} -m(k,\jj)+(s+ \rho)-m(k_0,\jj_0),0)}.
\end{equation*}
We finally get \eqref{eq:rx-R-th-R0-x} with Lemma \ref{lem:R-Rinf}, \eqref{eq:rxd-R-xd} and \eqref{eq:R-th-R0-derivee}. For $\nu \in \{0,\dots,k\}$, we have
\begin{eqnarray*}
\lefteqn{\nr{\pppg{rx}^\d \Rc_{\nu,\ll_\nu}(ir,z) \D_{G_0} \Rc_{k-\nu,\ll_\nu^*}^\infty(ir,z) \th_\sigma(z) \Rc^0_{\k_0,\jj_0}(ir,z) \pppg x^{-\d}}} \\
&& \lesssim \sum_{\abs{\a_1},\abs{\a_2} = 1} \nr{\pppg{rx}^\d  \Rc_{\nu,\ll_\nu}(ir,z) D^{\a_1} \pppg {x}^{-\d}} \nr{ D^{\a_2} \Rc_{k-\nu,\ll_\nu^*}^\infty(ir,z) \th_\sigma(z) \Rc^0_{\k_0,\jj_0}(ir,z) \pppg x^{-\d}} \\
&& \lesssim r^{\min(s-m(\nu,\ll_\nu)+1,0)} r^{\min(s + 1 - m(k-\nu,\ll_\nu^*) + (\s+\rho) -m(k_0,\jj_0),0)}.
\end{eqnarray*}
Estimate \eqref{eq:rx-R-th-R0-x} follows with \eqref{eq:min} applied with $\eta_0 = s$.
\end{proof}

To finish the proof of Proposition \ref{prop:regularity-Rz-A}, we have to replace $\pppg{rx}^\d$ by $\pppg{\A}^\d$ in \eqref{eq:x-R-rx}-\eqref{eq:rx-R-th-R0-x}. For this, we use again the elliptic regularity to compensate for the derivatives that appear in $\pppg{\A}^\d$.

\dimm{Ici à cause des $\gamma$ ça ne passe pas en l'état en dimension 2}

\begin{proposition} \label{prop:xTA}
Let $\d \geq 0$. Let $k \in \N$ and $\jj \in \{0,1,2\}^k$. Assume that $k \geq \d$. Let $\ell_1,\ell_2 \in \{0,1,2\}$. Then there exists $C > 0$ such that for $z \in \DD_+$ and $r=\abs z$, we have
\begin{equation*} 
\nr{\pppg {rx}^{-\d} \g_{\ell_1}(z) \Rc_{k,\jj}(ir,z) \gamma_{\ell_2}(z) \pppg {\A}^\d}_{\Lc(L^2)}\leq C  r ^{\ell_1 -m(k,\jj)+ \ell_2},
\end{equation*}
and
\begin{equation*}
\nr{\pppg {\A}^\d \g_{\ell_1}(z)   \Rc_{k,\jj}(ir,z) \g_{\ell_2}(z)  \pppg {rx}^{-\d}}_{\Lc(L^2)} \leq C  r ^{\ell_1 -m(k,\jj)+ \ell_2}.
\end{equation*}
\end{proposition}

Note that the assumption $k \geq \d$ is not optimal, but we have as many resolvents as we wish in Proposition \ref{prop:regularity-Rz-A}.

\begin{proof}
We start by proving by induction on $\k \in \N$ that, for any $
k \geq \k$ and $\mu \in \N^d$, 
\begin{equation} \label{eq:x-res-A}
\nr{\pppg {rx}^{-\k} \ad_{rx}^\m \big(\g_{\ell_1}(z) \Rc_{k,\jj}(ir,z) \gamma_{\ell_2}(z) \big)  \A^\k}_{\Lc(L^2)} \lesssim   r ^{\ell_1 -m(k,\jj)+ \ell_2}.
\end{equation}
With $\mu = 0$, we can deduce the first estimate of the proposition when $\d$ is an even integer. The general case will then follow by interpolation.
Estimate \eqref{eq:x-res-A} with $\k = 0$ is established by Propositions \ref{prop:Pa-Po-comm} and \ref{prop:reg-T} (notice that $m'(k,\jj) \geq 2 \geq \min(1,\ell_1) + \min(1,\ell_2)$).

Now let $\k \in \N^*$ and $\mu \in \N^d$, and assume that $k \geq \k$. We can write the commmutator $\ad_{rx}^\m \big( \g_{\ell_1}(z) \Rc_{k,\jj}(ir,z) \gamma_{\ell_2}(z)  \big)$ as a sum of terms of the form $$\ad_{rx}^{\m_1} \big( \g_{\ell_1}(z) \Rc_{k-1,\jj_1}(ir,z) \gamma_{j_k}(z)\big)\ad_{rx}^{\m_2}\big( R(ir) \gamma_{\ell_2}(z) \big)$$ where
\begin{equation}
    \label{cond:nu}
      \jj_1 = \{j_1,\dots,j_{k-1}\} \quad \text{ and } \quad \mu_1 + \mu_2 = \mu.
\end{equation}
Then we can write
\[\pppg {rx}^{-\k}  \ad_{rx}^{\m_1} \big( \g_{\ell_1}(z) \Rc_{k-1,\jj_1}(ir,z) \gamma_{j_k}(z)\big)\ad_{rx}^{\m_2}\big( R(ir) \gamma_{\ell_2}(z) \big) \A^\k\]
as a sum of terms of the form
\[
\pppg {rx}^{-\k} \ad_{rx}^{\m_1} \big( \g_{\ell_1}(z) \Rc_{k-1,\jj_1}(ir,z) \gamma_{j_k}(z)\big)  \A^p \ad_{\A}^{\k-p}\ad_{rx}^{\m_2}\big( R(ir) \gamma_{\ell_2}(z)  \big),
\]
where  $p \in \Ii 0 {\k}$. For $p\leq \k-1$, we apply the induction assumption to write
\begin{eqnarray*}
\lefteqn{\nr{\pppg {rx}^{-\k} \ad_{rx}^{\m_1} \big( \g_{\ell_1}(z) \Rc_{k-1,\jj_1}(ir,z) \gamma_{j_k}(z)\big)  \A^p  }_{\Lc(L^2)}} \\
&&\lesssim   \nr{\pppg {rx}^{-p} \ad_{rx}^{\m_1} \big( \g_{\ell_1}(z) \Rc_{k-1,\jj_1}(ir,z) \gamma_{j_k}(z)\big)  \A^p }_{\Lc(L^2)}\\
&& \lesssim r ^{\ell_1 -m(k-1,\jj_1)+ j_k}.
\end{eqnarray*}
By Propositions \ref{prop:Pa-Po-comm} and \ref{prop:reg-T}, we have
\[\nr{ \ad_{\A}^{\k-j}\ad_{rx}^{\m_2}\big( R(ir) \gamma_{\ell_2}(z)  \big)}_{\Lc(L^2)}\lesssim r ^{\ell_2 -2} .\]
Using the last two estimates, along with \eqref{cond:nu}, we obtain the estimate for $p \leq \k-1 $. Now we consider the term corresponding to $p = \k$. We have
\[
\A^\k = -\frac {id}{2} \A^{\k-1}  + \A^{\k-1} \sum_{\ell=1}^d r x_\ell \cdot r \inv D_\ell.
\]
The contribution of the first term is estimated as before. Now let $\ell \in \Ii 1 d$. By Propositions \ref{prop:Pa-Po-comm} and \ref{prop:reg-T} again, we have
\[ \nr{ r\inv D_\ell \ad_{rx}^{\m_2}\big( R(ir) \gamma_{\ell_2}(z)  \big)}_{\Lc(L^2)}\lesssim   r ^{\ell_2 -2}. \]
On the other hand, we have  by \eqref{eq:comm-Ak-xl}
\begin{eqnarray*}
\lefteqn{\pppg {rx}^{-\k}  \ad_{rx}^{\m_1} \big( \g_{\ell_1}(z) \Rc_{k-1,\jj_1}(ir,z) \gamma_{j_k}(z)\big) \A^{\k-1} rx_\ell}\\
&& = \pppg {rx}^{-\k} \ad_{rx}^{\m_1} \big( \g_{\ell_1}(z) \Rc_{k-1,\jj_1}(ir,z) \gamma_{j_k}(z)\big) rx_\ell (\A-i)^{\k-1} \\
&& = rx_\ell \pppg {rx}^{-\k}  \ad_{rx}^{\m_1} \big( \g_{\ell_1}(z) \Rc_{k-1,\jj_1}(ir,z) \gamma_{j_k}(z)\big)(\A-i)^{\k-1}\\
&& + \pppg {rx}^{-\k} \ad_{rx_\ell} \big( \ad_{rx}^{\m_1} \big( \g_{\ell_1}(z) \Rc_{k-1,\jj_1}(ir,z) \gamma_{j_k}(z)\big)\big)  (\A-i)^{\k-1}.
\end{eqnarray*}
Both terms are estimated with the induction assumption, thus proving \eqref{eq:x-res-A} and hence the first estimate of the proposition. The proof of the second estimate is similar.
\end{proof}

Now we can conclude the proof of Proposition \ref{prop:regularity-Rz-A}.

\begin{proof}[Proof of Proposition \ref{prop:regularity-Rz-A}]
We combine Propositions \ref{prop:elliptic-2} and \ref{prop:xTA}. We consider \eqref{eq:xRA-Schro-1}. Since $N$ can be chosen large we have, for some $k_1,k_2 \in \N$, $\jj_1 \in \{0,1,2\}^{k_1}$, $\jj_2 \in \{0,1,2\}^{k_2}$ and $\tilde \ell \in \{0,1,2\}$,
\[
\Rc_{N, \jj}(ir,z) = \Rc_{k_1,\jj_1}(ir,z) \g_{\tilde \ell}(z) \Rc_{k_2,\jj_2}(ir,z),
\]
with 
\[
m(k_1,\jj_1) - \tilde \ell + m(k_2,\jj_2) = m(N,\jj), \quad m(k_1,\jj_1) \geq s, \quad k_2 \geq \d.
\]
Then
\begin{eqnarray*}
\lefteqn{\big\| \pppg x^{-\d}  \Rc_{N, \jj}(ir,z)\gamma_{\ell}(z) \pppg {\A}^\d \big\|_{\Lc(L^2)}}\\
&& \leq \big\| \pppg x^{-\d}  \Rc_{k_1, \jj_1}(ir,z) \pppg {rx}^\d \big\|_{\Lc(L^2)}\big\| \pppg {rx}^{-\d}  \g_{\tilde \ell}(z) \Rc_{k_2, \jj_2}(ir,z)\gamma_{\ell}(z) \pppg {\A}^\d \big\|_{\Lc(L^2)} \\
&& \lesssim r^{s-m(k_1,\jj_1)} r^{\tilde \ell - m(k_2,\jj_2) + \ell} \lesssim r^{s- m(N,\jj) + \ell}.
\end{eqnarray*}
The other estimates are proved similarly.
\end{proof}

\section{Mourre commutators method} \label{sec:Mourre}

We set
\[
\DR^\pm = \set{z \in \DD_+ \, : \, \pm 2 \Re(z) \geq \abs z^2}, \quad \DR = \DR^+ \cup \DR^-.
\]

In this section, we prove Proposition \ref{prop:Mourre} about terms involving resolvents $R(z)$ with $z$ in $\DR$. For $z \in \DI$, Proposition \ref{prop:Mourre} holds with $\delta = 0$. This is a consequence of elliptic regularity, and more precisely of Proposition \ref{prop:regularity-Rz} applied without weights.

The proof of Proposition \ref{prop:Mourre} relies on the Mourre commutators method. We use the version given in \cite[Sec.5]{Royer24}, which we recall now.\\

We consider two Hilbert spaces $\Hc$ and $\Kc$ such that $\Kc$ is densely and continuously embedded in $\Hc$. We identify $\Hc$ with its dual, so that $\Kc \subset \Hc \subset \Kc^*$. Compared to \cite{Royer24}, we choose the convention that the dual of a Hilbert space is the set of continuous antilinear forms. Then the identification $\Hc \simeq \Hc^*$ and the natural embedding $\Ic : \Kc \to \Kc^*$ are linear.

Let $Q \in \Lc(\Kc,\Kc^*)$. We assume that $Q$ is boundedly invertible and that its imaginary part has a sign. For instance, we assume that it is non-positive, and we set
\begin{equation} \label{eq:Q+dissipative}
Q_+ := -\Im(Q)  = - \frac {Q-Q^*} {2i} \geq 0.
\end{equation}
In general, we say that an operator in $\Lc(\Kc,\Kc^*)$ is dissipative if its imaginary part is non-positive.\\

Let $\Asf$ be a selfadjoint operator on $\Hc$, with domain $\Dc_\Hc$. It can be seen as a bounded operator $\Asf_\Hc$ from $\Dc_\Hc$ (endowed with the graph norm) to $\Hc$. By restriction, $\Asf_\Hc$ defines a bounded operator $\Asf_\Kc$ from $\Dc_\Kc = \set{\f \in \Kc \cap \Dc_\Hc \, : \, \Asf_\Hc\f \in \Kc}$ to $\Kc$. Then $\Asf_\Kc^*$ defines a bounded operator $\Asf_{\Kc^*}$ from $\Dc_{\Kc^*} = \set {\f \in \Kc^* \, : \, \Asf_\Kc^* \f \in \Kc^*}$ to $\Kc^*$.

Let $\Kc_1,\Kc_2 \in \{\Kc,\Hc,\Kc^*\}$. We set $\Cc^0_\Asf(\Kc_1,\Kc_2) = \Lc(\Kc_1,\Kc_2)$ and for $S \in \Lc(\Kc_1,\Kc_2)$ we set $\ad_\Asf^0(S) = S$. Then, by induction on $n \in \N^*$, we say that $S \in \Cc^{n}_\Asf(\Kc_1,\Kc_2)$ if $S \in \Cc^{n-1}_\Asf(\Kc_1,\Kc_2)$ and the commutator $\ad_\Asf^{n-1}(S)  \Asf_{\Kc_1} - \Asf_{\Kc_2^*}^* \ad_\Asf^{n-1}(S)$, which defines at least a bounded operator from $\Dc_{\Kc_1}$ to $\Dc_{\Kc_2^*}^*$, extends to an operator $\ad_\Asf^n(S)$ in $\Lc(\Kc_1,\Kc_2)$. Then, we set
\[
\nr{S}_{\Cc^n_\Asf(\Kc_1,\Kc_2)} = \sum_{k=0}^n \big\| \ad_\Asf^k (S)\big\|_{\Lc(\Kc_1,\Kc_2)}.
\]
We write $\Cc^n(\Kc_1)$ for $\Cc^n(\Kc_1,\Kc_1)$. \\

Now we can recall the definition of a conjugate operator to $Q$. Compared to the version given in \cite{Royer24}, we add a parameter $\h$. We do not use it here (we choose $\h = 1$ in Proposition \ref{prop:Az-conjugate} below), but it is useful when dealing with high frequencies (see \cite{royer-mourre,BoucletRoy14}).

\begin{definition} \label{def-conjugate-operator}
Let $N \in \N^*$, $\h \in ]0,1]$ and $\Upsilon \geq 1$. We say that $\Asf$ is $(\h,\Upsilon)$-conjugate to $Q$ up to order $N$ if the following conditions are satisfied.
\begin{enumerate}[\rm (H1)]
\item \label{hyp-M-Kc} For $\f \in \Kc$ we have $\nr{\f}_\Hc \leq \Upsilon \nr{\f}_\Kc$.
\item \label{hyp-M-propagator} For all $\th \in [-1,1]$ the propagator $e^{-i\th \Asf} \in \Lc(\Hc)$ defines by restriction a bounded operator on $\Kc$.

\item \label{hyp-M-B} $Q$ belongs to $\Cc^{N+1}_\Asf(\Kc,\Kc^*)$ and
\begin{enumerate}
 \item $\nr{Q}_{\Lc(\Kc,\Kc^*)} \leq \Upsilon$,
 \item $\nr{\ad_\Asf(Q)}_{\Lc(\Kc,\Kc^*)} \leq \sqrt \h \Upsilon$,
 \item $\nr{\ad_\Asf^k(Q)}_{\Lc(\Kc,\Kc^*)} \leq \h \Upsilon$, for all $k \in \Ii 2 {N+1}$.
\end{enumerate}

\item \label{hyp-M-Qbot-Qpp}
There exist $Q_\bot \in  \Lc(\Kc,\Kc^*)$ dissipative, $\Qpp \in  \Lc(\Kc,\Kc^*)$ non-negative and $\Pi \in \Cc^1_\Asf(\Hc,\Kc)$ such that, with $\Pi_\bot = \Id_\Kc - \Pi \in \Lc(\Kc)$,
\begin{enumerate}
\item $Q = Q_\bot -i \Qpp$,
\item $\nr{\Qpp}_{ \Lc(\Kc,\Kc^*)} \leq \Upsilon$, $\nr{\Pi}_{\Cc^1_\Asf(\Hc,\Kc)} \leq \Upsilon$, and for $\f \in \Hc$ we have $\nr{\Pi \f}_\Kc \leq \Upsilon \nr{\Pi \f}_{\Hc}$,
\item $Q_\bot$ has an inverse $R_\bot \in  \Lc(\Kc^*,\Kc)$ which satisfies $\nr{\Pi_\bot R_\bot}_{ \Lc(\Kc^*,\Kc)} \leq \Upsilon$ and $\nr{R_\bot \Pi_\bot^*}_{ \Lc(\Kc^*,\Kc)} \leq \Upsilon$.
\end{enumerate}

\item \label{hyp-M-comm-pos}
There exists $\b \in [0,\Upsilon]$ such that if we set
\[
M =  i \ad_{\Asf}(Q) + \b \Qp \in  \Lc(\Kc,\Kc^*),
\]
then
\[
\nr{\ad_\Asf (\Pi^* M \Pi)}_{\Lc(\Kc,\Kc^*)} \leq \h \Upsilon,
\]
and, in the sense of quadratic forms on $\Hc$, we have
\begin{equation*}
\Pi^* \Re (M) \Pi \geq \frac \h \Upsilon \Pi^* \Ic \Pi.
\end{equation*}
\end{enumerate}
\end{definition}

We recall the estimates of $Q\inv$ given by the commutators method. See Theorem 5.2 in \cite{Royer24}. For the proof we follow \cite{Royer24} and check the dependence in $\h$ of all the estimates (see also \cite{royer-mourre}).

\begin{theorem} \label{th:mourre-simple}
Let $Q \in \Lc(\Kc,\Kc^*)$ be dissipative and boundedly invertible. Let $\Asf$ be a selfadjoint operator on $\Hc$. Let $N \in \N^*$, $\h \in ]0,1]$ and $\Upsilon \geq 1$. Assume that $\Asf$ is $(\h,\Upsilon)$-conjugate to $Q$ up to order $N$.
\begin{enumerate} [\rm (i)]
\item
Let $\d > \frac 12$. There exists $c > 0$, which only depends on $\Upsilon$ and $\d$, such that
\begin{equation} \label{eq:mourre-simple-1}
\nr{\pppg \Asf^{-\d} Q\inv  \pppg \Asf^{-\d}} _{\Lc(\Hc)}  \leq \frac c \h.
\end{equation}
\item Assume that $N \geq 2$ and let $\d_1,\d_2 \geq 0$ be such that $\d_1 + \d_2< N-1$. There exists $c > 0$, which only depends on $N$, $\Upsilon$, $\d_1$ and $\d_2$, such that
\begin{equation} \label{eq:mourre-simple-2}
\nr{\pppg{\Asf}^{\d_1} \1_{\R_-} (\Asf) Q\inv \1_{\R_+}(\Asf) \pppg{\Asf}^{\d_2}}_{\Lc(\Hc)} \leq \frac c \h.
\end{equation}
\item
Assume that $N \geq 2$ and let $\d \in \left] \frac 12 , N \right[$. There exists $c > 0$, which only depends on $N$, $\Upsilon$ and $\d$, such that
\begin{equation} \label{eq:mourre-simple-3}
\nr{\pppg{\Asf}^{-\d} Q\inv  \1_{\R_+}(\Asf) \pppg{\Asf}^{\d-1}}_{\Lc(\Hc)}  \leq \frac c \h
\end{equation}
and
\begin{equation} \label{eq:mourre-simple-4}
\nr{ \pppg{\Asf} ^{\d-1}\1_{\R_-}(\Asf) Q\inv  \pppg{\Asf}^{-\d}}_{\Lc(\Hc)}  \leq \frac c \h.
\end{equation}
\end{enumerate}
\end{theorem}

Using \cite[Lemma 5.5]{Royer24} we can convert the simple resolvent estimates into multiple resolvent estimates, with inserted factors. This gives the following result.

\begin{theorem} \label{th:mourre-multiple}
Let $\Asf$ be a selfadjoint operator on $\Hc$. Let $n \in \N^*$ and $N > n$. Let $Q_1,\dots,Q_n \in \Lc(\Kc,\Kc^*)$ be dissipative and boundedly invertible. We assume that for all $j \in \{1,\dots,n\}$ the selfadjoint operator $\Asf$ is $(\h_j,\Upsilon_j)$-conjugate to $Q_j$ up to order $N$ for some $\Upsilon_j \geq 0$ and $h_j \in ]0,1]$. Let $B_0,\dots,B_{n} \in C^N_\Asf(\Hc)$. Let
\[
\mathcal T = B_0 Q_1\inv B_1 Q_2\inv \dots B_{n-1} Q_n\inv B_n.
\]
Let $\d \in \left] n-\frac 12 , N \right[$ and let $\d_1,\d_2 \geq 0$ be such that $\d_1 + \d_2< N-n$. There exists $C > 0$ which only depends on $N$, $\Upsilon_1,\dots,\Upsilon_n$, $\d$, $\d_1$ and $\d_2$ such that
\begin{equation} \label{eq:mourre-multiple-1}
\nr{\pppg \Asf^{-\d}  \mathcal T   \pppg \Asf^{-\d}} _{\Lc(\Hc)}  \leq \frac {C}{\h_1 \dots \h_n} \nr{B_0}_{C_\Asf^N(\Hc)} \dots \nr{B_n}_{C_\Asf^N(\Hc)},
\end{equation}
\begin{equation} \label{eq:mourre-multiple-2}
\nr{\pppg{\Asf}^{\d_1} \1_{\R_-} (\Asf) \, \mathcal T \, \1_{\R_+}(\Asf) \pppg{\Asf}^{\d_2}}_{\Lc(\Hc)} \leq \frac {C}{\h_1 \dots \h_n}\nr{B_0}_{C_\Asf^N(\Hc)} \dots \nr{B_n}_{C_\Asf^N(\Hc)},
\end{equation}
\begin{equation} \label{eq:mourre-multiple-3}
\nr{\pppg{\Asf}^{-\d} \, \mathcal T \,  \1_{\R_+}(\Asf) \pppg{\Asf}^{\d-1}}_{\Lc(\Hc)}  \leq \frac {C}{\h_1 \dots \h_n}\nr{B_0}_{C_\Asf^N(\Hc)} \dots \nr{B_n}_{C_\Asf^N(\Hc)},
\end{equation}
\begin{equation} \label{eq:mourre-multiple-4}
\nr{ \pppg{\Asf} ^{\d-1}\1_{\R_-}(\Asf) \, \mathcal T \, \pppg{\Asf}^{-\d}}_{\Lc(\Hc)}  \leq \frac {C}{\h_1 \dots \h_n}\nr{B_0}_{C_\Asf^N(\Hc)} \dots \nr{B_n}_{C_\Asf^N(\Hc)}.
\end{equation}
\end{theorem}

\begin{proof}
We apply \cite[Lemma 5.5]{Royer24} with $(2n+1)$ factors $T_1,\dots,T_{2n+1}$. More precisely, we set
\begin{itemize}
\item $T_{2j+1} = B_j / \nr{B_j}_{C_\Asf^N(\Hc)}$, $\nu_{2j+1} = \s_{2j+1} = 0$, for all $j \in \Ii 0 n$,
\item $T_{2j} = h_j Q_j\inv$, $\nu_{2j} = 1$, $\s_{2j} = \d - n + 1$ for all $j \in \Ii 1 n$,
\item $\Th_j = \pppg \Asf$, for all $j \in \Ii 0 {2n+1}$,
\item $\Pi_j^+ = \1_{\R_+}(\Asf)$, $\Pi_j^- = \1_{\R_-^*}(\Asf)$, for all $j \in \Ii 0 {2n+1}$.
\end{itemize}
The estimates for $Q_j\inv$ are given by Theorem \ref{th:mourre-simple}. For an operator $B_j$, we apply \cite[Proposition 3.11]{Royer24}. Then Lemma \cite[Lemma 5.5]{Royer24} provides a collection of constants $\Cc = \{ C ; (C_{\d_-,\d_+}) ; (C_\d^-) ; (C_\d^+) \}$ which only depend on $N$, $\Upsilon_1,\dots,\Upsilon_n$, $\d$, $\d_1$ and $\d_2$ such that \eqref{eq:mourre-multiple-1}-\eqref{eq:mourre-multiple-4} hold. The conclusions follow.
\end{proof}

For the proof of Proposition \ref{prop:Mourre}, we apply this abstract result in $\Hc = L^2$ to the resolvent $R(z)$ with the generator of dilations \eqref{def:A} as the conjugate operator. For each $z \in \DR$, we choose $\Kc = H^1_z$. Everything depends on $z$, but the key in the following statement is that $\Upsilon$ does not.

\begin{proposition}\label{prop:Az-conjugate}
Let $N \in \N$. There exist $\Upsilon \geq 1$ such that for all $z \in \DR^+$ the generator of dilations $\A$ is $(1,\Upsilon)$-conjugate to $\abs z^{-2} P(z) \in \Lc(H^1_z,H_z^{-1})$ up to order $N$.
\end{proposition}

Compared to the proof of Proposition 5.4 in \cite{Royer24}, we add the contribution of the damping term, and we check that we can work with the usual generator of dilations $\A$.

For the proof, we will use the Helffer-Sj\"ostrand formula, which we recall now. Let $H$ be a selfadjoint operator on a Hilbert space $\Hc$, $m \geq 2$ and let $\vf \in C_0^\infty(\R)$. Let $\psi \in C_0^\infty(\R,[0,1])$ be supported on $[-2,2]$ and equal to 1 on $[-1,1]$. The almost analytic extension $\tilde \vf$ of $\vf$ is defined by
\[
\tilde \vf(\tau+i\mu) = \psi\left( \frac \mu {\pppg \tau} \right) \sum_{k=0}^m \phi^{(k)}(\tau) \frac {(i\mu)^k}{k!}.
\]
Then we have
\begin{equation*}
\vf(H) = - \frac 1 \pi \int_{\C} \frac {\partial \tilde \phi}{\partial \bar \z} (\z) (H-\z)\inv \, d\l(\z),
\end{equation*}
where $\l$ is the Lebesgue measure on $\C$. See for instance \cite[Section 8]{dimassis} or \cite{davies95}.

\begin{proof} For $z \in \DR^+$ we set $Q_z = \abs z^{-2} P(z)$.

\stepp The assumption (H\ref{hyp-M-Kc}) holds for any $\Upsilon \geq 1$ since $\nr{\cdot}_{L^2} \leq \nr{\cdot}_{H^1_z}$ for any $z \in \DR^+$. For (H\ref{hyp-M-propagator}), we recall that for $\th \in [-1,1]$, $u \in L^2(\R^d)$ and $x \in \R^d$, we have $(e^{i\th A} u)(x) = e^{\frac {d\th}2} u ( e^{\th} x)$. Thus, if $u$ belongs to $H_z^1 = H^1(\R)$, then so does $e^{i\th A} u$, and
$\| e^{i\th A} u \|_{H^1_z(\R^d)} \leq e \nr u_{H^1_z}$,
which gives (H\ref{hyp-M-propagator}).
By Proposition \ref{prop:Pa-Po-comm}, we have $Q_z \in C^N_z(H^1_z,H\inv_z)$ with norm independent of $z \in \DR^+$, and (H\ref{hyp-M-B}) follows.

\stepp We set
\[
\PR(z)= \Re(P(z)) = \abs z^2 \Re(Q_z) = -\D_G +aw\Im(z)- w \Re(z^2).
\]
This defines a selfadjoint operator on $L^2$ with domain $H^2$, and also a bounded operator from $H_z^1$ to $H_z\inv$. We have
\begin{align*}
[\PR(z),iA]
& = -2\D_G + K_1(z),
\end{align*}
where
\begin{align*}
K_1(z) = \sum_{j,k} \partial_j ((x\cdot \nabla) G_{j,k}(x)) \partial_k - \Im(z) (x\cdot \nabla) (aw) + \Re(z^2) (x \cdot \nabla) w.
\end{align*}
Moreover,
\[
-2 \D_G = 2 \PR(z) + 2 \Re(z^2) + K_2(z),
\quad 
K_2(z) = - 2 \Im(z) aw + 2 \Re(z^2) (w-1).
\]
Setting $K(z) = K_1(z) + K_2(z)$, we obtain
\begin{equation} \label{eq:commutator}
[\PR(z),iA] = 2 \PR(z) + 2 \Re(z^2) + K(z).
\end{equation}
We fix $\rho \in ]0,\rho_0[$. Then there exists $C_1 > 0$ such that 
\begin{equation} \label{eq:estim-Kz}
\nr{ K(z) \pppg {zx}^{\rho}}_{\Lc(H_z^1,H_z\inv)} \leq C_1 \abs z^2.
\end{equation}

\stepp Now we construct the operator $\Pi_z$ which appears in (H\ref{hyp-M-Qbot-Qpp}) and (H\ref{hyp-M-comm-pos}).
Let $\vf \in C_0^\infty(\R,[0,1])$ be equal to 1 on $[ -1,1]$ and supported in $]-2,2[$. For $\eta \in ]0,1]$ and $\l \in \R$, we set $\phi_\eta (\l) = \phi(\l / \eta)$. Then for $z \in \DR^+$, we set
\[
\Pi_{\eta,z} = \vf_\eta \big( \Re(Q_z) \big) = \vf_\eta \left( \frac {\PR(z)}{\abs z^2} \right).
\]
We also set 
\[
\Pi_{\eta,z}^0 =\vf_\eta \big( \Re(Q_z^0) \big), \quad  Q_z^0 = \frac {-\D - z^2}{\abs z^2}.
\]
For all $\eta \in ]0,1/32]$ we have
\begin{equation} \label{eq:PiPRPi}
2 \nr{\Pi_{2\eta,z} \PR(z) \Pi_{2\eta,z}}_{\Lc(L^2)} \leq \frac {\abs z^2} 4.
\end{equation}
 Moreover, with the Helffer-Sjöstrand formula, we can check that there exists $C_2 > 0$ such that for all $z \in \DR^+$ and $\eta \in ]0,1/32]$, we have
\begin{equation} \label{eq:norm-Piz}
\nr{\Pi_{2\eta,z}}_{\Lc(H_z\inv,H_z^1)} \leq C_2.
\end{equation}

\stepp We prove that there exists a family $(\eta_z)_{z \in \DR^+}$ in $]0,1/32]$ such that $\eta_* := \inf \eta_z > 0$ and, for all $z \in \DR^+$,
\begin{equation} \label{eq:Pi-K-Pi}
\nr{\Pi_{2\eta_z,z} K(z) \Pi_{2\eta_z,z}}_{\Lc(L^2)} \leq \frac {\abs z^2}4.
\end{equation}
Let $z_0 \in \overline{\DR^+}$. We prove that there exist $\eta_{z_0} \in ]0,1/32]$ and a neighborhood $\Vc_{z_0}$ of $z_0$ in $\C$ such that for all $z \in \DR^+ \cap \Vc_{z_0}$ we have \eqref{eq:Pi-K-Pi} with $\eta_z = \eta_{z_0}$.

First assume that $z_0 \neq 0$. Since $\Pi_{1,z_0} \in \Lc(H\inv,L^2)$, $K(z_0) \pppg x^{\rho} \pppg D\inv \in \Lc(L^2,H\inv)$, $\pppg D \pppg x^{-\rho} \pppg D^{-2}$ is a compact operator on $L^2$ and $\pppg D^2 \Pi_{1,z_0} \in \Lc(L^2)$, the operator
\[
\Pi_{1,z_0} K(z_0) \Pi_{1,z_0} = \Pi_{1,z_0} \big(K(z_0) \pppg x^{\rho} \pppg D\inv \big) \big( \pppg D \pppg x^{-\rho} \pppg D^{-2} \big) \big( \pppg D^2 \Pi_{1,z_0} \big)
\]
is compact on $L^2$. Since 0 is not an eigenvalue of $\PR(z_0)$, the operator $\Pi_{\eta,z_0}$ goes weakly to 0 as $\eta$ goes to 0, so there exists $\eta_{z_0} \in ]0,1/32]$ such that
\[
\nr{\Pi_{2\eta_{z_0},z_0} K(z_0) \Pi_{2\eta_{z_0},z_0}}_{\Lc(L^2)} \leq \frac {\abs{z_0}^2} 8.
\]
By continuity with respect to $z$, there exists a neighborhood $\Vc_{z_0}$ of $z_0$ such that, for all $z \in \DR^+ \cap \Vc_{z_0}$,
\begin{equation*}
\nr{\Pi_{2\eta_{z_0},z} K(z) \Pi_{2\eta_{z_0},z}}_{\Lc(L^2)} \leq \frac {\abs z^2}4.
\end{equation*}

\stepp Now, we proceed similarly for $z$ in a neighborhood of 0. Let $\tau_0 \in \big[\frac 1 {\sqrt 2} , 1 \big]$. Since $\pppg {D} \pppg{x}^{- \rho} \pppg {D}^{-2}$ is compact as an operator from $L^2$ to $H^1$ and $\pppg D^2 \vf_{\eta_0}(\frac {-\D-\tau_0^2}{16})$ goes weakly to 0 as $\eta_0$ goes to 0, there exists $\eta_0 \in \big]0,\frac 1 8 \big]$ such that
\begin{equation}
\begin{aligned}
\nr{\pppg {D_z} \pppg{zx}^{-\rho} \vf_{\eta_0} \left(\frac {-\D-\tau_0^2 \abs z^2}{16 \abs z^2} \right)}_{\Lc(L^2)}
& = \nr{\pppg D \pppg{x}^{-\rho} \vf_{\eta_0} \left(\frac {-\D-\tau_0^2}{16}\right)}_{\Lc(L^2)}\\
& \leq \frac 1 {8 C_1 C_2},
\end{aligned}
\end{equation}
where $C_1$ and $C_2$ are the constants which appear in \eqref{eq:estim-Kz} and \eqref{eq:norm-Piz}. If $\abs{\Re(z^2/\abs z^2) - \tau_0^2} \leq 8 \eta_0$ we have
\[
\Pi_{2\eta_0,z}^0 =  \vf_{\eta_0} \left( \frac {-\D-\tau_0^2 \abs z^2}{16 \abs {z}^2} \right) \Pi_{2\eta_0,z}^0,
\]
so
\begin{equation} \label{eq:phi-phi}
\nr{\pppg{zx}^{-\rho } \Pi_{2\eta_0,z}^0}_{\Lc(L^2,H^1_z)} \leq \frac 1 {8 C_1 C_2}.
\end{equation}
Since $\big[\frac 1 {\sqrt 2} , 1 \big]$ is compact, we can choose $\eta_0$ so small that \eqref{eq:phi-phi} holds for any $z \in \DR^+$.

By the Helffer-Sj\"ostrand formula (applied with $m \geq 3$) and the resolvent identity, the difference $\Pi_{2\eta_0,z} - \Pi_{2\eta_0,z}^0$ can be rewritten as
\[
\frac {1} {\pi} \int_{\C} \frac {\partial \tilde \vf_{2\eta_0}}{\partial \bar \z}(\zeta) \left( \Re(Q_z) - \z \right)\inv \Re(Q_z-Q_z^0)   \left( \Re(Q_z^0) - \z \right)\inv \, d\l(\z).
\]
We can check that for $z \in \DD_+$ and $\z \in 5\DD \setminus \R_+$, we have
\begin{equation} \label{eq:res-PR}
\nr{\left(\Re(Q_z)- \z \right)\inv}_{\Lc(H_z^{-1},H_z^{1})} + \nr{\left(  \Re(Q_z^0) -  \z \right)\inv}_{\Lc(L^2,H^2_z)} \lesssim \frac{1} {\abs{\Im(\z)}}.
\end{equation}
On the other hand,
\[
\nr{\Re(Q_z-Q_z^0)}_{\Lc(H_z^{1+\rho},H_z^{-1})} \lesssim \abs z^{\rho}.
\]
This proves
\[
 \nr{\left( \Re(Q_z) - \z \right)\inv \Re(Q_z-Q_z^0)   \left( \Re(Q_z^0) - \z \right)\inv }_{\Lc(L^2,H_z^{1})} \lesssim \frac {\abs{z}^{\rho}} {\abs{\Im(\z)}^2} .
\]
Since $\partial_{\bar \z} \tilde \vf_{2\eta_0}$ is compactly supported and decays faster than $\abs{\Im(\z)}^2$ near the real axis, we deduce
\begin{equation} \label{eq:Piz-Piz0}
\nr{\Pi_{2\eta_0,z} - \Pi_{2\eta_0,z}^0}_{\Lc(L^2,H_z^{1})} \lesssim \abs{z}^{\rho}.
\end{equation}
Then, by \eqref{eq:phi-phi} and \eqref{eq:Piz-Piz0}, there exists $r_0 \in ]0,1]$ such that for $z \in \DR$ with $\abs z \leq r_0$ we have
\begin{equation} \label{eq:zx-Pi}
\nr{\pppg{zx}^{-\rho} \Pi_{2\eta_0,z}}_{\Lc(L^2,H^1_z)} \leq \frac 1 {4 C_1 C_2}.
\end{equation}
With \eqref{eq:estim-Kz} and \eqref{eq:norm-Piz} we get
\begin{equation} \label{eq:K-Piz}
\nr{\Pi_{2\eta_0,z} K(z)\Pi_{2\eta_0,z}}_{\Lc(L^2)} \leq \frac {\abs z^2}{4}.
\end{equation}
This concludes the proof of \eqref{eq:Pi-K-Pi}. The fact that we can choose the family $(\eta_z)$ bounded away from 0 comes from the fact that $\eta_z$ can be chosen constant in a neighborhood of any $z_0$ in the compact $\overline{\DR^+}$.

\stepp
Now we can prove (H\ref{hyp-M-Qbot-Qpp}) and (H\ref{hyp-M-comm-pos}). We begin with (H\ref{hyp-M-comm-pos}). We choose $\b=0$. Let $z \in \DR^+$. By \eqref{eq:commutator}, \eqref{eq:PiPRPi} and \eqref{eq:Pi-K-Pi}, we have for all $z \in \DR^+$
\[
\Pi_{2\eta_z,z} [\PR(z),iA] \Pi_{2\eta_z,z} \geq 2 \Re(z^2) \Pi_{2\eta_z,z}^2 - \frac {\abs z^2} 2.
\]
Since $2\Re(z^2) \geq \abs z^2$ we get after composition by $\Pi_{\eta_z,z}$ on both sides
\[
\Pi_{\eta_z,z} [\PR(z),iA] \Pi_{\eta_z,z} \geq \frac {\abs z^2} 2 \Pi_{\eta_z,z}^2.
\]
This gives (H\ref{hyp-M-comm-pos}) with $\Pi_z := \Pi_{\eta_z,z}$.

\stepp
By the Helffer-Sj\"ostrand formula and Lemma \ref{lem:commutators} we have
\begin{align*}
\nr{[\Pi_z,iA_z]}_{\Lc(H\inv_z,H_z^1)}
& \lesssim \int_\C \abs{\frac {\partial \tilde \vf_{\eta_z}}{\partial \bar \z}(\z)}  \nr{\left[\left(\Re(Q_z) -\z \right)\inv,iA \right]}_{\Lc(H\inv_z,H_z^1)} \, d\l(\z) \\
& \lesssim  \nr{ \big[\Re(Q_z),iA \big]}_{\Lc(H^1_z,H_z\inv)}\\
& \lesssim 1.
\end{align*}
Now we set
\[
Q_\bot(z) = \frac {\PR(z) - i \Im(z^2) w_{\min}}{\abs z^2} \in \Lc(H_z^1,H_z\inv),
\]
where $w_{\min} = \min_{x \in \R} w(x) > 0$. Then
\begin{align*}
\Qpp(z)
 := i \big(Q_z - Q_\bot(z) \big)
 =\abs z^{-2} \big[\Im(z^2) (w-w_{\min})+\Re(z) aw\big]
\end{align*}
is non-negative, $Q_\bot(z)$ is boundedly invertible and by the functional calculus we have
\[
\nr{(1-\Pi_z) Q_\bot(z) \inv}_{\Lc(H_z\inv,H^1_z)} = \nr{Q_\bot(z) \inv(1-\Pi_z) }_{\Lc(H_z\inv,H^1_z)} \leq \frac 1 {\eta_*}.
\]
This gives (H\ref{hyp-M-Qbot-Qpp}) and the proof is complete.
\end{proof}

Now we can prove Proposition \ref{prop:Mourre}.

\begin{proof}[Proof of Proposition \ref{prop:Mourre}]
We begin with \eqref{estim-Mourre-wave-1}.
If $k=0$, it is given by Theorem \ref{th:mourre-simple} applied with $\Asf=\A$ and $Q=\frac{1}{|z|^2} P(z)$.
Assume that $k=1$, so that $\Rc_{k,\jj} = R(z) \g_{j_1}(z) R(z)$ for some $j_1 \in\{0,1\}$. If $j_1=0$ we simply apply Theorem \ref{th:mourre-multiple} with $Q_1= Q_2 = \frac{1}{|z|^2} P(z)$ and $B_1 = \g_0(z)$. When $j_1 = 1$, we want to see $\g_1(z)$ as an operator from $H_z^{s}$ to $H_z^{s-1}$ for some $s$, since then it is of size $O(z)$. Then we cannot apply directly Theorem \ref{th:mourre-multiple}, where the inserted factors are bounded operators on $L^2$. Thus, we apply once more the resolvent identity  \eqref{eq:res-identity} to get elliptic regularity. Setting $r = \abs z$ and $\tilde \g_1(z) = R(ir) \g_1(z) R(ir)$, we can write
\begin{eqnarray*}
\lefteqn{R(z) \gamma_{1}(z) R(z)
=   \big[1+R(z)\gamma_2(z) \big] \tilde \gamma_{1}(z) \big[ 1+ \gamma_2(z) R(z)\big]}\\
&&= \tilde \g_1(z) + \tilde \g_1(z) \gamma_2(z) R(z)+ R(z)\gamma_2(z) \tilde \g_1(z)+ R(z)\gamma_2(z)   \tilde \g_1(z) \gamma_2(z) R(z).
\end{eqnarray*}
By Propositions \ref{prop:Pa-Po-comm} and \ref{prop:Rari-H1}, we have
\[
\nr{\tilde \g_1(z)}_{\Cc_{\A}^N(L^2)} \lesssim \frac 1 {\abs z^3}, \quad \nr{\tilde \g_1(z) \gamma_2(z)}_{\Cc_{\A}^N(L^2)} \lesssim \frac 1 {\abs z},\quad \nr{ \gamma_2(z)\tilde \g_1(z)}_{\Cc_{\A}^N(L^2)} \lesssim \frac 1 {\abs z},
\]
and
\[
\nr{\gamma_2(z)   \tilde \g_1(z) \gamma_2(z)}_{\Cc_{\A}^N(L^2)} \lesssim \abs z.
\]
\dimm{On utilise $d\geq 3$ encore ici, et c'est peut-être l'endroit le plus problématique}
Then, we can apply Theorem \ref{th:mourre-multiple} for each term and deduce \eqref{estim-Mourre-wave-1} for $k=1$ and $z \in \DR^+$. We proceed similarly for $k \geq 2$, the details are omitted.

We turn to \eqref{estim-Mourre-wave-2}. We use the same trick to see the inserted factor $\th_\s(z)$ as an operator of size $O(\abs z^{\s+\rho})$. Assume for simplicity that $k = k_0 = 0$. Setting $\tilde \th_{\sigma }(z)= R(ir)\th_\sigma(z)  R_0(ir)$, we can write, with the resolvent identities \eqref{eq:res-identity} and  \eqref{eq:res0-identity},
\begin{eqnarray*}
\lefteqn{R(z)\th_\sigma(z) R_0(z)= \big[ 1+ R(z)\gamma_2(z) \big] \tilde\th_{\sigma}(z)  \big[ 1+ \gamma^0_2(z) R_0(z)\big]}\\
&&= \tilde\th_{\sigma}(z)+ \tilde\th_{\sigma}(z) \gamma^0_2(z) R_0(z)+ R(z)\gamma_2(z)  \tilde\th_{\sigma}(z)+ R(z)\gamma_2(z)\tilde\th_{\sigma}(z)\gamma^0_2(z) R_0(z).
\end{eqnarray*}
By Propositions \ref{prop:Pa-Po-comm}, \ref{prop:Rari-H1} and Proposition \ref{prop:reg-Ra} for $R_0(z)$, we have
\begin{eqnarray*}
\lefteqn{\| \gamma_2(z)\tilde\th_{\sigma}(z)\gamma_2^0(z) \|_{\Cc_{\A}^N(L^2)}
\lesssim \abs z^4 \|\tilde\th_{\sigma}(z)\|_{\Cc_{\A}^N(L^2,H^1_z)} }\\
&& \lesssim \abs z^4 \nr{R(ir)}_{\Cc_{\A}^N(H_z\inv,H_z^1)} \nr{\th_\sigma(z)}_{\Cc_{\A}^N(H^{1+\rho}_z,H_z\inv)} \nr{R_0(ir)}_{\Cc_{\A}^N(L^2,H^{1+\rho}_z)}\\
&& \lesssim \abs z^{\s+\rho} .
\end{eqnarray*}
Similarly, 
\[
\| \tilde \th_\s(z) \|_{\Cc_{\A}^N(L^2)} \lesssim \abs z^{(\s+\rho) - 4},
\]
\[
\| \tilde\th_{\sigma}(z)\gamma_2^0(z) \|_{\Cc_{\A}^N(L^2)} \lesssim \abs z^{(\s+\rho) - 2} , \quad
\| \gamma_2(z)\tilde\th_{\sigma}(z) \|_{\Cc_{\A}^N(L^2)}\lesssim \abs z^{(\s+\rho) - 2}.
\]
We deduce \eqref{estim-Mourre-wave-2} for $z \in \DR^+$, and the proof is then similar when $k,k_0 \geq 1$.

Then Proposition \ref{prop:Mourre} is proved for $z \in \DR^+$. By a duality argument, we deduce \eqref{estim-Mourre-wave-1} and \eqref{estim-Mourre-wave-2} for $z \in \DR^-$, and the proof is complete.
\end{proof}

\appendix

\section{Additional proof} \label{sec:additional-proofs}

In this appendix we prove Lemma \ref{lem:T-N-resolventes}.

\begin{proof}[Proof of Lemma \ref{lem:T-N-resolventes}]
The result is proved by induction on $N \in \N$. We begin with the case $N = 0$. If $k = 0$ then the claim is directly given by the resolvent identity \eqref{eq:res-identity}. Now assume that $k \geq 1$. We have
\begin{align*}
\Rc_{k,\jj}(z)=R(z)\gamma_{j_1}(z) R(z)... \gamma_{j_k}(z) R(z)  =R(z)\gamma_{j_1}(z)\Rc_{k-1,\jj'}(z).
\end{align*}
where we have set $\jj' = (j_2,\dots,j_k) \in \{0,1\}^{k-1}$.
By \eqref{eq:res-identity} we get
\begin{align*}
\Rc_{k,\jj}(z)
 = R(ir)  \gamma_{j_1}(z) \Rc_{k-1,\jj_1}(z) + R(ir) \gamma_2(z)  \Rc_{k,\jj}(z).
\end{align*}
Both terms are of the form \eqref{eq:T-N-res-2}. Now, let $N \geq 1$ and assume that the result is proved for $N-1$. Then it is enough to prove that an operator of the form \eqref{eq:T-N-res-1} or \eqref{eq:T-N-res-2} (with $N-1$ instead of $N$) can be written as a sum of terms of the form \eqref{eq:T-N-res-1} or \eqref{eq:T-N-res-2}. A term of the form \eqref{eq:T-N-res-1} is already in suitable form. We consider a term of the form \eqref{eq:T-N-res-2} (for $N-1$).
Let $\jj_1 = (j_{1,1},\dots,j_{1,N-1}) \in \{0,1,2\}^{N-1}$, $\jj_2= (j_{2,1},\dots,j_{2,k_2}) \in \{0,1\}^{k_2}$ (with $k_2 \leq k$) and $\ell \in \{0,1,2\}$.
If $k_2 = 0$, we have by \eqref{eq:res-identity}
\begin{eqnarray*}
\lefteqn{\Rc_{N-1,\jj_1}(ir,z) \gamma_{\ell}(z) R(z)}\\
&& = \Rc_{N-1,\jj_1}(ir,z) \gamma_{\ell}(z) R(ir)
+ \Rc_{N-1,\jj_1}(ir,z) \gamma_{\ell}(z) R(ir) \gamma_2(z) R(z) \\
&& = \Rc_{N,\jj_3}(ir,z) + \Rc_{N,\jj_3}(ir,z) \gamma_2(z) R(z),
\end{eqnarray*}
where $\jj_3 = (j_{1,1},\dots,j_{1,N-1},\ell) \in \{0,1,2\}^{N}$.
If $k_2 \geq 1$, we similarly have with $\jj_4= (j_{2,2},\dots,j_{2,k_2}) \in \{0,1\}^{k_2-1}$
\begin{eqnarray*}
\lefteqn{\Rc_{N-1,\jj_1}(ir,z) \gamma_{\ell}(z) \Rc_{k_2,\jj_2}(z)}\\
&&=\Rc_{N-1,\jj_1}(ir,z) \gamma_{\ell}(z) R(z) \gamma_{j_{2,1}}(z) \Rc_{k_2-1,\jj_4}(z)\\
&& = \Rc_{N,\jj_3}(ir,z)  \gamma_{j_{2,1}}(z) \Rc_{k_2-1,\jj_4}(z) + \Rc_{N,\jj_3}(ir,z)  \gamma_2(z) R(z) \gamma_{j_1}(z)  \Rc_{k_2-1,\jj_4}(z).\\
&& = \Rc_{N,\jj_3}(ir,z)  \gamma_{j_{2,1}}(z) \Rc_{k_2-1,\jj_4}(z) + \Rc_{N,\jj_3}(ir,z)  \gamma_2(z) \Rc_{k_2,\jj_2}(z).
\end{eqnarray*}
In both cases, we get two terms of suitable form.
\end{proof}

\subsection*{Acknowledgments}
This work has been supported by the Labex CIMI, Toulouse, France, under grant ANR-11-LABX-0040-CIMI.

\end{document}